\documentclass{article}

\usepackage{amsmath}
\usepackage{amssymb}
\usepackage{amsthm}
\usepackage{boldtensors}
\usepackage{geometry}
\usepackage{graphicx}
\usepackage{lineno}
\usepackage{microtype}
\usepackage{multirow}
\usepackage{url} %this package should fix any errors with URLs in refs.
\usepackage{xcolor}
\theoremstyle{plain}
\newtheorem{remark}{Remark}
\usepackage[boxed]{algorithm}
\usepackage{algpseudocode}
\usepackage{booktabs}
\usepackage{makecell}
\usepackage{tikz}

\usepackage[backend=bibtex,maxcitenames=2]{biblatex}
\usepackage[font=footnotesize]{caption}
\usepackage{subcaption}

\newcommand{\numDOFs}{N}
\newcommand{\numModes}{M}
\newcommand{\numSnaps}{\ensuremath{N_t}}
\newcommand{\numSamps}{\ensuremath{N_s}}

\newcommand{\numOverlap}{\ensuremath{N_o}}
\newcommand{\numBound}{\ensuremath{N_b}}
\newcommand{\timeIdx}{n}
\newcommand{\iterIdx}{k}
\newcommand{\varIdx}{v}

\newcommand{\ode}[2]{\ensuremath{\frac{\text{d} #1}{\text{d} #2}}}
\newcommand{\pde}[2]{\ensuremath{\frac{\partial #1}{\partial #2}}}
\newcommand{\pdeTwo}[2]{\ensuremath{\frac{\partial^2 #1}{\partial {#2}^2}}}
\newcommand{\solVar}{u}
\newcommand{\solVec}{\ensuremath{\mathbf{\solVar}}}
\newcommand{\solMat}{\ensuremath{\mathbf{\MakeUppercase{\solVar}}}}
\newcommand{\forceVar}{f}
\newcommand{\forceVec}{\ensuremath{\mathbf{\forceVar}}}
\newcommand{\nonlinOpVar}{n}
\newcommand{\nonlinOp}{\ensuremath{\mathcal{\MakeUppercase{\nonlinOpVar}}}}
\newcommand{\nonlinOpVec}{\ensuremath{\mathbf{\MakeUppercase{\nonlinOpVar}}}}
\newcommand{\resVar}{r}
\newcommand{\resVec}{\ensuremath{\mathbf{\resVar}}}
\newcommand{\resFunc}[1]{\ensuremath{\resVec \left( #1 \right)}}

\newcommand{\solVecRom}{\ensuremath{\mathbf{\widetilde{\solVar}}}}
\newcommand{\solVecCoef}{\ensuremath{\mathbf{\widehat{\solVar}}}}
\newcommand{\solVecCent}{\ensuremath{\mathbf{\overline{\solVar}}}}
\newcommand{\dummyVec}{\ensuremath{\mathbf{y}}}
\newcommand{\resOp}{\ensuremath{\mathbf{A}}}
\newcommand{\jacobMat}{\ensuremath{\mathbf{J}}}
\newcommand{\sampMat}{\ensuremath{\mathbf{S}}}
\newcommand{\standVec}{\ensuremath{\mathbf{e}}}

\newcommand{\trialVar}{\phi}
\newcommand{\trialVarCap}{\Phi}
\newcommand{\trialVec}{\ensuremath{\boldsymbol{\trialVar}}}
\newcommand{\trialBasis}{\ensuremath{\mathbf{\trialVarCap}}}

\newcommand{\testVarCap}{\Psi}

\newcommand{\testBasis}{\ensuremath{\mathbf{\testVarCap}}}
\newcommand{\projVarCap}{\Pi}
\newcommand{\projMat}{\ensuremath{\mathbf{\projVarCap}}}
\newcommand{\dummyBasis}{\ensuremath{\mathbf{\widetilde{X}}}}
\newcommand{\leftSVecMat}{\ensuremath{\mathbf{X}}}
\newcommand{\rightSVecMat}{\ensuremath{\mathbf{Y}}}
\newcommand{\EValMat}{\ensuremath{\mathbf{\Sigma}}}

\newcommand{\inROne}[1]{\ensuremath{\in \mathbb{R}^{#1}}}
\newcommand{\inRTwo}[2]{\ensuremath{\in \mathbb{R}^{#1 \times #2}}}
\newcommand{\funcMap}[3]{\ensuremath{#1 : #2 \rightarrow #3}}
\newcommand{\defEq}{\ensuremath{:=}}

\begin{document}

\title{The role of interface boundary conditions and sampling strategies for Schwarz-based coupling of projection-based reduced order models}

\author{Christopher R. Wentland$^1$, Francesco Rizzi$^{2}$,  Joshua Barnett$^{3}$,  Irina Tezaur$^1$\thanks{Email: ikalash@sandia.gov}
  \\
  \\
  \small $^1$Sandia National Laboratories, Livermore, CA, USA\\
  \small $^2$NexGen Analytics, Sheridan, WY USA\\
   \small $^3$Cadence Design Systems, San Jose, CA, USA\\
}

\date{\today}

\maketitle

% Abstract
\begin{center}
{\bf Abstract}
\end{center}

%\begin{center}
%    {\bf Abstract} \\
%\end{center}

This paper presents and evaluates a framework for the coupling of subdomain-local projection-based reduced order models (PROMs) using the Schwarz alternating method following a domain decomposition (DD) of the spatial domain on which a given problem of interest is posed. In this approach, the solution on the full domain is obtained via an iterative process in which a sequence of subdomain-local problems are solved, with information propagating between subdomains through transmission boundary conditions (BCs). We explore several new directions involving the Schwarz alternating method aimed at maximizing the method's efficiency and flexibility, and demonstrate it on three challenging two-dimensional nonlinear hyperbolic problems: the shallow water equations, Burgers' equation, and the compressible Euler equations. We demonstrate that, for a cell-centered finite volume discretization and a non-overlapping DD, it is possible to obtain a stable and accurate coupled model utilizing Dirichlet-Dirichlet (rather than Robin-Robin or alternating Dirichlet-Neumann) transmission BCs on the subdomain boundaries. We additionally explore the impact of boundary sampling when utilizing the Schwarz alternating method to couple subdomain-local hyper-reduced PROMs. Our numerical results suggest that the proposed methodology has the potential to improve PROM accuracy by enabling the spatial localization of these models via domain decomposition, and achieve up to two orders of magnitude speedup over equivalent coupled full order model solutions and moderate speedups over analogous monolithic solutions. 

% Introduction

\section{Introduction} \label{sec:intro}

% \ikt{Possible reviewers: Ivan Prusak, Ramon Codina, Giulia Sambotaro, Michel Bergmann, Alejandro Diaz, Matthias Heinkerschloss, John Jakeman?, Pete Bosler?}

Despite improved algorithms and high performance computing (HPC) architectures, the high computational costs associated with performing high-fidelity modeling and simulation (mod/sim) of complex physical or engineered systems can preclude analyses. Recent years have seen efforts to mitigate this problem through the development of data-driven models, which can run in a fraction of the time of high-fidelity simulations. Projection-based model order reduction (PMOR) is a promising approach for constructing data-driven models that have a low online computational cost and that, importantly, maintain the underlying partial differential equations (PDEs) being modeled. PMOR works to reduce the computational complexity of numerical simulations by restricting the search of the solution to a low-dimensional space spanned by a reduced basis constructed from a limited number of high-fidelity simulations and/or physical experiments/observations. The corresponding projection-based reduced order model (PROM) is obtained by performing a Galerkin or Petrov--Galerkin~\cite{Carlberg:2017} projection of the governing PDEs onto a reduced subspace, and approximating the state variables in a low-dimensional subspace spanned by a reduced basis. Popular approaches for computing linear (or affine) subspaces include proper orthogonal decomposition (POD)~\cite{Holmes:1996, Sirovich:1987}, balanced POD~\cite{Rowley:2005, Willcox:2002} and the reduced basis method~\cite{Rozza:2011, Veroy:2005}. Recently, nonlinear manifold approaches have emerged and have been shown to deliver the same accuracy as their traditional affine counterparts but using a much smaller dimension, especially for challenging convection-dominated problems having a slowly decaying Kolmogorov $n$-width~\cite{Pinkus:1985}. Examples of nonlinear manifold bases include autoencoders~\cite{Lee:2020, Kim:2022}, quadratic approximation manifolds~\cite{Barnett:2022, Geelen:2023} or artificial neural network- (ANN-)augmented bases~\cite{Barnett:2023}.

While recent years have seen numerous advancements in the development of PROMs and other data-driven models, the use of these models in predictive settings raises some fundamental questions about their numerical properties. Unfortunately, traditional PROMs are known to suffer from a lack of robustness, stability, and accuracy, especially in the predictive regime. Moreover, because traditional PROMs are usually characterized by spatially-global trial spaces, they typically require relatively high-dimensional representations to preserve model accuracy for complex problems (e.g., problems with a slowly decaying Kolmogorov $n$-width~\cite{Pinkus:1985}), and cannot be easily adjusted to accommodate geometry modifications. While combining PROMs with high-fidelity models can potentially mitigate these challenges, a unified and rigorous theory for integrating these models in a ``plug-and-play'' fashion into existing mod/sim toolchains is lacking at the present time.

In this paper, we develop a new methodology for improving the predictive power of PROMs and facilitating the rigorous incorporation of these models into mod/sim workflows using coupling and domain decomposition. Our approach is based on three key ingredients: (i) the decomposition of the underlying spatial geometry into overlapping or non-overlapping subdomains, (ii) the creation of subdomain-local PROMs that can better capture difficult solution features than global PROMs, and (iii) the rigorous coupling of subdomain-local PROMs to each via the Schwarz alternating method~\cite{Schwarz:1870}. The Schwarz alternating method is an iterative domain decomposition-based approach based on a very simple idea: if the solution to a PDE is known in two or more regularly-shaped subdomains comprising a more complex domain, these local solutions can be used to iteratively build a solution on the more complex domain, with information propagating between subdomains through carefully-specified transmission boundary conditions (BCs). The approach presented in this paper builds on our past work involving the development of the Schwarz alternating method as a means to achieve concurrent multi-scale coupling of high-fidelity solid mechanics models~\cite{Mota:2017, Mota:2022}, and as a novel contact enforcement algorithm~\cite{Mota:2023}. It was demonstrated in these works that the approach has a number of advantages: (i) it is minimally intrusive to implement in existing HPC software frameworks; (ii) it is capable of coupling regions with different mesh resolutions, different element types, and different time integration schemes (e.g., implicit and explicit), all without introducing any artifacts exhibited by alternative coupling methods; and (iii) it possesses rigorous convergence properties/guarantees~\cite{Mota:2017, Mota:2022}. We develop the Schwarz alternating method as a means to enable the coupling of subdomain-local PROMs, and study the method's numerical properties on several hyperbolic conservation laws. While the method can additionally couple subdomain-local PROMs with subdomain-local full order models (FOMs), that specific use case is not considered here for the sake of brevity. We focus our attention on PROMs constructed using the POD/least-squares Petrov--Galerkin (LSPG) method~\cite{Carlberg:2011}, as POD/LSPG is generally the method of choice for hyperbolic problems, like those considered here. To enable computational efficiency for PROMs of nonlinear systems, we employ hyper-reduction via collocation on randomized sample meshes. We describe the specific contributions and differentiating features of our approach after surveying the literature for related past work.

\subsection{Related past work} \label{sec:past_work}

The notion of coupling subdomain-local PROMs with each other and/or with FOMs has been gaining attention in the PMOR literature in recent years. Existing coupling approaches can be classified as falling into two categories: (i) monolithic coupling methods, and (ii) iterative coupling methods. Methods from category (i) are typically more intrusive to implement but can possess more favorable convergence properties and runtimes than methods falling into the category (ii).

The vast majority of monolithic coupling methods are based on formulations that first define appropriate interface compatibility constraints, and then enforce these constraints using different techniques. Among the earliest approaches for coupling subdomain-local PROMs with each other and with FOMs using Lagrange multipliers is the work of Lucia \textit{et al.}~\cite{Lucia:2003} (for POD-based PROMs) and Maday \textit{et al.}~\cite{Maday:2004} (for PROMs constructed via the reduced basis method). Both~\cite{Lucia:2003} and~\cite{Maday:2004} impose compatibility using a mortar-type method that ``glues'' together non-overlapping subdomains using Lagrange multipliers approximated using low-order polynomials. Several subsequent works approach the coupling problem as a means to ``stitch'' together a variety of composable PROM ``tiles.'' In~\cite{Wicke:2009}, continuity between tiles is enforced by duplicating the degrees of freedom on the interfaces and constraining their normal components to be equal. A Lagrange multiplier-free approach that is conceptually similar to~\cite{Wicke:2009} is the reduced basis method with domain decomposition and finite elements (RDF)~\cite{Iapichino:2016}. In this non-overlapping approach, traditional finite element basis functions are defined on the boundaries between PROM domains so as to ensure continuity and provide a sort of finite element enrichment, with the goal of gluing together networks of repetitive blocks. Another recent work that accomplishes monolithic PROM-FOM coupling without Lagrange multipliers is~\cite{Baiges:2013}. Here, a discontinuous Galerkin (DG) formulation is assumed and coupling is achieved through a careful definition of numerical fluxes at discrete cell boundaries. In~\cite{Hoang:2021}, Hoang \textit{et al.} develop an algebraically non-overlapping method for coupling subdomain-local LSPG PROMs with each other. Compatibility conditions at subdomain boundaries can be imposed either strongly or weakly using Lagrange multipliers, and subdomain training can be performed independently. In a similar vein, Chung \textit{et al.}~\cite{Chung:2024} recently proposed their so-called ``data-driven finite element method'' (DD-FEM), which trains spatially-local reduced order models (ROMs) corresponding to small-scale unit components, tiles a computational domain with these components, and employs ideas from DG methods to glue the components together. In this method, the interface constraints are only enforced weakly through interface operators introduced to penalize constraint violations, thereby eliminating the need to introduce Lagrange multipliers. In~\cite{Diaz:2024}, Diaz \textit{et al.} extends the approach in~\cite{Hoang:2021} to nonlinear manifold autoencoder-based ROMs, which can be more efficient for problems in the slowly-decaying Kolmogorov $n$-width regime. In~\cite{deCastro:2022, deCastro:2023}, de Castro \textit{et al.} develop a Lagrange multiplier-based explicit synchronous partitioned scheme for performing PROM-PROM and PROM-FOM coupling. One of the advantages of this method is that, when combined with explicit time-integration, the subdomain-local problems can be solved independently through a clever utilization of the Schur complement. Moreover, when utilizing independently-constructed PROM bases for the interfacial and interface degrees of freedom, it is shown that the approach leads to a provably nonsingular Schur complement independent of the underlying mesh size and/or reduced basis dimension~\cite{deCastro:2023}.

The Schwarz alternating method falls into a category of iterative coupling methods. As mentioned above, these methods are usually less intrusive to implement in HPC codes, but may be less efficient than monolithic methods. Among the earliest Schwarz-based DD approaches for PROM-PROM and PROM-FOM coupling focuses on Galerkin-free PROMs. In~\cite{Buffoni:2007}, Buffoni \textit{et al.} propose three methods to perform Galerkin-free PROM-FOM and PROM-PROM coupling: (i) a Schur iteration where the solution of the PROM is obtained by a projection step in the space spanned by the POD modes, (ii) a Dirichlet-Dirichlet iteration in the frame of a classical Schwarz method, and (iii) an approach obtained by minimizing the residual norm of the canonical approximation in the space spanned by the POD modes. Galerkin-free ROM-FOM and ROM-ROM couplings are also explored by Cinquegrana \textit{et al.}~\cite{Cinquegrana:2011} and Bergmann \textit{et al.}~\cite{Bergmann:2018}. The former approach considers overlapping DD in the context of a Schwarz-like iteration scheme, but, unlike our approach, requires matching meshes at the subdomain interfaces~\cite{Cinquegrana:2011}. The latter approach, known as zonal Galerkin-free POD, defines an optimization problem which minimizes the difference between the POD reconstruction and its corresponding FOM solution in the overlapping region between a ROM and a FOM domain~\cite{Bergmann:2018}. A true POD-greedy/Galerkin non-overlapping Schwarz method for the coupling of PROMs developed for the specific case of symmetric elliptic PDEs is presented by Maier \textit{et al.} in~\cite{Maier:2014}. In the recent work~\cite{Iollo:2022}, Iollo \textit{et al.} demonstrated that the classical overlapping Schwarz iteration can be interpreted as an optimization-based coupling, which can be solved directly in the method known as ``one-shot Schwarz.'' Optimization-based domain decomposition frameworks for coupling subdomain-local PROMs/FOMs have also been considered in several other recent works, e.g.,~\cite{Prusak:2023} and~\cite{Hawkins:2024}.

Some emerging technologies in PMOR coupling are applicable to both monolithic and iterative coupling formulations. The approaches summarized above assume a predefined domain decomposition (DD) and PROM/FOM subdomain assignment, which is consistent with the scope of this paper. It is worth noting, however, that PROM-PROM and PROM-FOM coupling methods with on-the-fly basis and/or DD adaptation are also possible and beginning to emerge in the PMOR literature. In~\cite{Corigliano:2013, Corigliano:2015}, Corigliano \textit{et al.} develop an adaptive non-overlapping Lagrange multiplier-based coupling method for nonlinear elasto-plastic and multi-physics problems. Here, on-the-fly PROM adaptation and PROM/FOM switching is performed through a plastic check during the reduced analysis. Hybrid PROM–FOM coupling in the context of solid mechanics applications is considered also in~\cite{Kerfriden:2012, Kerfriden:2013, Radermacher:2014}. In~\cite{Kerfriden:2012, Kerfriden:2013}, a local/global model reduction strategy for the simulation of quasi-brittle fracture is developed, whereas, in~\cite{Radermacher:2014}, an adaptive sub-structuring (or DD) approach is proposed that not only enables on-the-fly adaptation of the PROM basis, but also on-the-fly substructuring/DD changes. Finally, in~\cite{Huang:2022}, Huang \textit{et al.} develop a component-based modeling framework that can flexibly integrate PROMs and FOMs for different components or domain decompositions, with an eye towards large-scale combustion problems. It is demonstrated that accuracy can be enhanced by incorporating basis adaptation ideas from~\cite{Peherstorfer:2020, Huang:2023}.

Researchers have additionally begun in recent years to integrate ideas from machine learning (ML) into PROM-PROM and PROM-FOM coupling formulations. In~\cite{Ahmed:2021}, Ahmed \textit{et al.} present a hybrid PROM-FOM approach in which a long short-term memory (LSTM) network is introduced at the subdomain interfaces and utilized to perform the multi-model coupling. In order to improve efficiency, Bochev \textit{et al.}~\cite{Bochev:2024} replace the Schur complement solve required in the explicit synchronous partitioned presented in~\cite{deCastro:2022, deCastro:2023} with a dynamic mode decomposition surrogate trained from interface flux data. The second direction involves the application to non-intrusive PROMs. In a similar vein to~\cite{Bochev:2024}, in~\cite{Discacciati:2024p2}, Discacciati \textit{et al.} propose a non-intrusive PMOR approach, which approximates the boundary response maps arising from a non-overlapping DD method using a combination of dimensionality reduction techniques and interpolation/regression approaches (e.g., kernel interpolation methods and ANNs). In~\cite{Farcas:2023}, Farcas \textit{et al.} develop a methodology for the coupling of nonintrusive and subdomain-local operator inference (OpInf) models, which can both improve accuracy and reduce memory/computation requirements. A similar approach for coupling subdomain-local OpInf models with each other and with FOMs was recently introduced in~\cite{Moore:2024}. In~\cite{LiD3M, LiDeepDDM}, the Schwarz alternating method was extended to enable the DD-based coupling of physics-informed neural networks. The methods proposed in these works, known as D3M~\cite{LiD3M} and DeepDDM~\cite{LiDeepDDM} respectively, inherit the benefits of DD-based ROM–ROM couplings, but are developed primarily for the purpose of improving the efficiency of the neural network training process and reducing the risk of over-fitting. A similar approach, but considering both weak and strong formulations of the Schwarz BCs, was explored by Snyder \textit{et al.} in~\cite{Snyder:2023}. The Schwarz alternating method has also been used for online coupling of independently pre-trained subdomain-localized neural network-based models, e.g. in~\cite{Wang:2022}, which develops a transferable framework for solving boundary value problems (BVPs) via deep neural networks that can be trained once and used forever for various unseen domains and BCs.

\subsection{Contributions and differentiating features} \label{sec:contribs}

While our manuscript builds on some of our earlier work in developing the Schwarz alternating method for FOM-FOM~\cite{Mota:2017, Mota:2022} and PROM-PROM/PROM-FOM~\cite{Barnett:2022} coupling in solid mechanics, it includes several novel contributions.

Whereas our past work focused on solid mechanics and finite element-based discretizations, here we consider challenging two-dimensional (2D) nonlinear hyperbolic conservation law problems, including problems with strong shocks, discretized using a cell-centered finite volume method (CCFV). While Dirichlet boundary conditions can be imposed strongly in FEM-based discretizations and PROMs based on top of these discretizations~\cite{Barnett:2022, Gunzburger:2007}, all boundary conditions in CCFV methods are enforced weakly using ghost cells, auxiliary computational cells that extend the physical domain for the purpose of boundary condition imposition. It is well-known that the non-overlapping Schwarz alternating method is in general non-convergent if Dirichlet-Dirichlet transmission boundary conditions are specified at the Schwarz boundaries; to achieve convergence, one must use either alternating Dirichlet-Neumann~\cite{LionsNonOverlap1988} or Robin-Robin~\cite{Zanolli:1987} transmission conditions. Implementing these conditions can be cumbersome in HPC codes, as the creation of carefully-constructed transmission operators is required~\cite{Mota:2023}. We demonstrate that the use of ghost cells for boundary condition imposition, which has the effect of introducing overlap into an otherwise non-overlapping domain decomposition, can improve the flexibility and simplify the implementation of the non-overlapping version of the Schwarz alternating method by enabling the use of Dirichlet-Dirichlet coupling conditions at the Schwarz boundaries. This is a novel result, to the best of our knowledge.

Our past work involving Schwarz for coupling has focused on the classical ``multiplicative'' variant of the method~\cite{Mota:2017, Mota:2022, Barnett:2022}. A well-known disadvantage of multiplicative Schwarz stems from the fact that the method solves subdomain-local problems in sequence until convergence is reached, which limits computational acceleration. This may be remedied by a variant of the method known as ``additive'' Schwarz, where subdomain-local problems are solved concurrently and boundary data are communicated asynchronously. Our numerical results demonstrate that, with subdomain-local PROMs coupled using additive Schwarz, it is possible to achieve up to two orders of magnitude speedup over equivalent coupled FOM solutions and moderate speedups over analogous monolithic FOM solutions. Our numerical experiments additionally reveal that a PROM-PROM coupling of a given size is typically more accurate than a monolithic PROM of comparable size defined on the full problem domain. This result suggests that the proposed coupling methodology has the potential to improve PROM accuracy by enabling the spatial localization of these models via domain decomposition.

Lastly, while our previous work on Schwarz for PROM-PROM coupling~\cite{Barnett:2022} considered nonlinear PDEs and thus included hyper-reduction to maintain efficiency, attention was restricted to one-dimensional problems, for which each subdomain had at most two Schwarz boundary nodes. In moving to the multi-dimensional case, where it is possible to have hundreds or thousands of Schwarz boundary mesh elements, it is important to develop a hyper-reduction strategy for sampling the Schwarz boundary in a way that maintains both accuracy and efficiency. We demonstrate that the optimal boundary sampling rate for applying Schwarz to hyperbolic conservation law PROMs with hyper-reduction is problem-specific, and, for the problems considered here, varies between sampling every two boundary cells (for problems exhibiting traveling shocks) and every ten boundary cells (for problems with smooth solutions).

The remainder of this paper is organized as follows. In Section~\ref{sec:schwarz}, we describe the Schwarz alternating method in the context of a general hyperbolic PDE, introducing several new ingredients relative to our past related work~\cite{Mota:2017, Mota:2022, Barnett:2022}, namely some details regarding additive Schwarz and the application of the method to CCFV discretizations. In Section~\ref{sec:rom}, we overview PMOR based on the POD/LSPG method with hyper-reduction, and present the details relevant to applying the Schwarz alternating method to PROMs constructed using this approach. We evaluate our method in the context of three 2D nonlinear hyperbolic conservation laws in Section~\ref{sec:results}: the shallow water equations, Burgers' equation, and the Euler equations. Conclusions are offered in Section~\ref{sec:conc}, which also includes a brief discussion of potential future work.

% Schwarz method
\section{The Schwarz alternating method for domain decomposition-based coupling} \label{sec:schwarz}

\begin{figure}
    \begin{minipage}{0.49\linewidth}
        \includegraphics[width=0.9\textwidth]{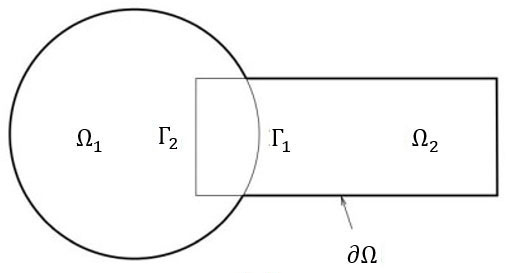}
        \subcaption{Overlapping}
        \label{subfig:dd_over}
    \end{minipage}
    \begin{minipage}{0.49\linewidth}
        \includegraphics[width=0.9\textwidth]{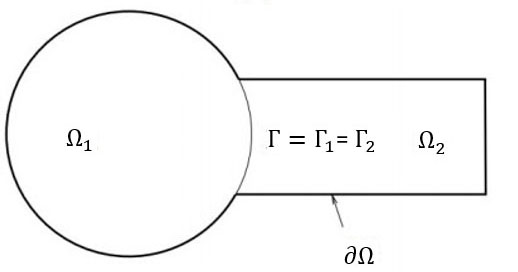}
        \subcaption{Non-overlapping}
        \label{subfig:dd_nonover}
    \end{minipage}
    \caption{Illustration showing overlapping and non-overlapping domain
    decomposition.}
    \label{fig:dd}
\end{figure}

In this section, we describe the Schwarz alternating method for domain decomposition-based coupling for a generic nonlinear PDE of the form 
\begin{equation} \label{eq:generic_pde}
\begin{array}{rcrl}
    \dot{\solVar}(x,t) + \nonlinOp(\solVar(x,t)) &=& \forceVar(x), & \text{in } \Omega \times [0, T], \\
    \solVar(x,t) &=& g(x), & \text{on } \partial \Omega \times [0,T]\\
    \solVar(x,0) & = & v(x), & \text{in } \Omega,
    \end{array}
\end{equation}
where $\Omega \inROne{d}$ is an open bounded domain for $d=1,2,3$ with boundary $\partial \Omega$, $\nonlinOp(\cdot)$ is a nonlinear operator, $\forceVar(x)$ and $g(x)$ are given source and boundary condition functions, respectively, $v(x)$ defines the initial condition, and $T > 0$. As we illustrate below, the Schwarz alternating method can be applied to both overlapping and non-overlapping domain decompositions, shown without loss of generality for the specific case of two subdomains in Figure~\ref{fig:dd}.

\subsection{The additive Schwarz iteration procedure} \label{sec:schwarz_procedure}

Consider the PDE~\eqref{eq:generic_pde} in $\Omega \times I_{\timeIdx}$ for some specified time-interval $I_{\timeIdx} = [t_{\timeIdx}, \ t_{\timeIdx+1}]$, where $T \ge t_{\timeIdx+1} > t_{\timeIdx} \ge 0$ with initial condition $\solVar(x, \ t_{\timeIdx}) = v(x)$ for a given function $v$. Assuming without loss of generality a two subdomain decomposition like those shown in Figure~\ref{fig:dd}, the additive Schwarz algorithm concurrently solves the following subdomain-local problems
\begin{equation} \label{eq:generic_schwarz_iter_Omega1}
    \left \{
    \begin{array}{rcll}
    \dot{\solVar}_1^{(\iterIdx)} + \nonlinOp(u_1^{(k)}) &=& f^{(k)}, & \text{in } \Omega_1 \times [t_{\timeIdx}, \ t_{\timeIdx+1}] \\ 
     \solVar_1^{(\iterIdx)} &=& g, & \text{on } (\partial \Omega \cap \overline{\Omega_1}) \times [t_{\timeIdx}, \ t_{\timeIdx+1}]\\
     \mathcal{B}_{11} \solVar_1^{(\iterIdx)} &=& \mathcal{B}_{12} \solVar_2^{(\iterIdx-1)},  &\text{on } \Gamma_1 \times [t_{\timeIdx}, \ t_{\timeIdx+1}], 
    \end{array}
    \right.
\end{equation}
and 
\begin{equation} \label{eq:generic_schwarz_iter_Omega2}
    \left \{
    \begin{array}{rcll}
        \dot{\solVar}_2^{(\iterIdx)} + \nonlinOp(\solVar_2^{(\iterIdx)}) &=& \forceVar^{(\iterIdx)}, & \text{in } \Omega_2 \times [t_{\timeIdx}, \ t_{\timeIdx+1}] \\ 
        \solVar_2^{(\iterIdx)} &=& g, & \text{on } (\partial \Omega \cap \overline{\Omega_2}) \times [t_{\timeIdx}, \ t_{\timeIdx+1}]\\
        \mathcal{B}_{22} \solVar_2^{(\iterIdx)} &=& \mathcal{B}_{21} \solVar_1^{(\iterIdx-1)},  &\text{on } \Gamma_2 \times [t_{\timeIdx}, \ t_{{\timeIdx}+1}], 
    \end{array}
    \right.
\end{equation}
for Schwarz iteration $\iterIdx = 1, 2, ...$ subject to the aforementioned initial conditions $\solVar_1(x,0) = v|_{\Omega_1}$ and $\solVar_2(x,0) = v|_{\Omega_2}$. Here, $\solVar_i$ for $i = 1, 2$ denotes the solution in subdomain $\Omega_i$, $\Gamma_i$ is the so-called Schwarz boundary (see Figure~\ref{fig:dd}), and $\mathcal{B}_{ij}$ for $i, j = 1, 2$ is a linear operator acting along the Schwarz boundary $\Gamma_i$ on the solution $\solVar_j$. It is common to set $\solVar_1^{(0)} = v|_{\partial \Omega_1}$ on $\Gamma_1$ and $\solVar_2^{(0)} = v|_{\partial \Omega_2}$ on $\Gamma_2$ when initializing the Schwarz iteration process, to ensure solution compatibility with the initial condition. The iterative process in~\eqref{eq:generic_schwarz_iter_Omega1} and~\eqref{eq:generic_schwarz_iter_Omega2} continues until a set of pre-determined criteria are met. In the present work, convergence criteria are based on the Euclidean norm of the solution differences between consecutive Schwarz iterations; that is, Schwarz is deemed converged when $\epsilon_{\text{abs}}^{(\iterIdx)} < \delta_{\text{abs}}$ and $\epsilon_{\text{rel}}^{(\iterIdx)} < \delta_{\text{rel}}$ for some pre-specified tolerances $\delta_{\text{abs}}, \delta_{\text{rel}} > 0$, where 
\begin{equation} \label{eq:conv_criterion_abs}
    \epsilon_{\text{abs}}^{(\iterIdx)} := \sqrt{|| \solVar_1^{(\iterIdx)} - \solVar_1^{(\iterIdx-1)}||^2_2 + || \solVar_2^{(\iterIdx)} - \solVar_2^{(\iterIdx-1)}||^2_2},
\end{equation}
and 
\begin{equation} \label{eq:conv_criterion_rel}
    \epsilon_{\text{rel}}^{(k)} := \sqrt{\frac{|| \solVar_1^{(\iterIdx)} - \solVar_1^{(\iterIdx-1)}||^2_2}{||\solVar_1^{(\iterIdx)}||^2_2} + \frac{|| \solVar_2^{(\iterIdx)} - \solVar_2^{(\iterIdx-1)}||^2_2}{||\solVar_2^{(\iterIdx)}||^2_2}},
\end{equation}
for Schwarz iteration $\iterIdx = 1, 2,...$.

Numerous past works have investigated the definitions of the operators $\mathcal{B}_{ij}$ to ensure well-posedness and convergence of each Schwarz subproblem and the Schwarz iteration process (the reader is referred to~\cite{Gander:2008} for a nice summary). Table~\ref{table:Schwarz_BCs} summarizes several choices of boundary operators defining a set of convergent Schwarz transmission conditions on the boundaries $\Gamma_i$ for the specific case of a two subdomain domain decomposition. Whereas Dirichlet-Dirichlet transmission conditions on $\Gamma_i$ will give rise to a convergent Schwarz iteration process provided the overlap region is non-empty~\cite{Lions:1988} ($\Omega_1 \cap \Omega_2 \neq \emptyset$; see Figure~\ref{subfig:dd_over}, non-overlapping domain decompositions ($\Omega_1 \cap \Omega_2 = \emptyset$; see Figure~\ref{subfig:dd_nonover} typically require alternating Dirichlet-Neumann~\cite{Zanolli:1987, Funaro:1988, Cote:2005, Kwok:2014} or Robin-Robin~\cite{LionsNonOverlap1988,Deng:2003, Lui:2001} transmission conditions on $\Gamma$. In Table~\ref{table:Schwarz_BCs}, projection operators $P_{ij}$ from subdomain $\Omega_j$ onto boundary $\Gamma_i$ are introduced to enable application of the method to the case where $\Omega_1$ and $\Omega_2$ are discretized using non-conformal meshes. As described in more detail in~\cite{Mota:2022, Mota:2023}, the projection operator from one subdomain to another consists of the simple application of the existing finite element interpolation functions. That is, any time the value of a field is needed at a point at the boundary in one subdomain, a search is performed to determine the element containing that point in the other subdomain, and the field value is calculated by using the nodal values of that element and the corresponding interpolation functions.

\begin{table}[ht!]
    \centering
    \caption{Summary of boundary operators defining convergent Schwarz transmission boundary conditions in~\eqref{eq:generic_schwarz_iter_Omega1} and~\eqref{eq:generic_schwarz_iter_Omega2} for a two subdomain domain decomposition.  $\partial_{n_i}$ denotes the normal derivative to boundary $\Gamma_i$, $P_{ij}$ denotes a projection operator from subdomain $\Omega_j$ onto boundary $\Gamma_i$, and $I$ denotes the identity.}
    \begin{tabular}{ccccc}
        \toprule
        DD Type & $\mathcal{B}_{11}$ & $\mathcal{B}_{12}$ & $\mathcal{B}_{21}$ & $\mathcal{B}_{22}$\\
        \midrule
        Overlapping  &  $I$ & $P_{12}$ & $I$ & $P_{21}$ \\
        \multirow{ 2}{*}{Non-overlapping} & $I$ & $P_{12}$ & $\partial_{n_2}$ & $\partial_{n_2}P_{21}$ \\
         & $\partial_{n_1} + I$ & $(\partial_{n_1} + I)P_{12}$  & $(\partial_{n_2} + I)P_{21}$ & $\partial_{n_2} + I$ \\
        \bottomrule
    \end{tabular}
    \label{table:Schwarz_BCs}
\end{table}

\begin{remark}
It is important to note that alternate transmission conditions besides the ``vanilla'' Dirichlet-Dirichlet, Robin-Robin and alternating Dirichlet-Neumann Schwarz transmission conditions described above are possible. In the so-called optimal Schwarz alternating method, relaxation is introduced into the transmission conditions~\cite{Gander:2008, Deng:2003, Lui:2001, Gerardo-Giorda:2013, Nataf:1997, Gander:2003, Halpern:2008}. The relaxation parameters are chosen so as to optimize the convergence of the Schwarz algorithm. It is also possible to define the $\mathcal{B}_{ij}$ in~\eqref{eq:generic_schwarz_iter_Omega1}--\eqref{eq:generic_schwarz_iter_Omega2} as non-local Dirichlet-to-Neumann (or Steklov-Poincare) operators~\cite{Gander:2008, Nataf:1994}. These advanced Schwarz formulations are beyond the scope of this paper and are not discussed in detail.
\end{remark} 

\begin{remark}
The interested reader may wonder when one may want to use an overlapping (Figure~\ref{subfig:dd_over} vs. a non-overlapping DD (Figure~\ref{subfig:dd_nonover}) when applying Schwarz. It has been demonstrated that the overlapping Schwarz variant tends to have more favorable convergence properties, meaning that the method typically converges in fewer Schwarz iterations~\cite{Gander:2008}. The non-overlapping Schwarz variant is more flexible and is the logical option for problems with a defined interface, e.g., multi-material or multi-physics (e.g., fluid-structure interaction) problems, but typically less efficient. Moreover, whereas Dirichlet boundary conditions are straightforward to enforce in most discretizations, a disadvantage of non-overlapping Schwarz is that the specification of Neumann or Robin BCs requires the design and implementation of effective transfer operators, which can be challenging to implement in HPC codes~\cite{Mota:2023}.
\end{remark} 

While the PDEs considered herein are dynamic, our discussion thus far has focused primarily on the spatial aspects of Schwarz-based coupling. For long-time dynamic problems defined on a large time-interval $[0, \ T]$, it is highly inefficient to integrate each subdomain $\Omega_i$ to time $T$ before synchronizing boundary conditions at the Schwarz boundary. If each subdomain is advanced forward in time using the same time-step $\Delta t$, a natural mitigation to this difficulty is to converge the Schwarz iteration process on a time-step by time-step basis, as illustrated in Algorithm~\ref{alg:schwarz}.

\begin{algorithm}[ht!]
\centering 
\caption{Additive Schwarz alternating method for a generic transient PDE~\eqref{eq:generic_pde} and a two subdomain domain decomposition. We are assuming that each subdomain is advanced forward in time using the same time-step.}

\begin{algorithmic}[1]
    \State $n \gets 0$
    \State $u_1^{(0)} = v|_{\Omega_1}$ and  $u_2^{(0)}= v|_{\Omega_2}$    \Comment initialize solutions in $\Omega_1$ and $\Omega_2$
    \Repeat
    \Comment{time-stepper}
        \State $\iterIdx \gets 1$
        \Comment initialize Schwarz iteration
        \Repeat
        \Comment{Schwarz loop}
            \State Compute $\mathcal{B}_{12}u_2^{(\iterIdx-1)}$
            \Comment interpolate $u_2^{(\iterIdx)}$ and its normal derivative to $\Gamma_1$
            \State Solve~\eqref{eq:generic_schwarz_iter_Omega1} in time interval $I_{\timeIdx}$ for $u_1^{(\iterIdx)}$
            \State Compute $\mathcal{B}_{21}u_1^{(\iterIdx-1)}$
            \Comment interpolate $u_1^{(\iterIdx)}$ and its normal derivative to $\Gamma_2$
            \State Solve~\eqref{eq:generic_schwarz_iter_Omega2} in time interval $I_{\timeIdx}$ for $u_2^{(\iterIdx)}$
            \State $\iterIdx \gets \iterIdx + 1$
            \Comment increment Schwarz iteration
        \Until $\epsilon_{\text{abs}}^{(\iterIdx)} < \delta_{\text{abs}}$ and $\epsilon_{\text{rel}}^{(\iterIdx)} < \delta_{\text{rel}}$
        \Comment{check convergence criteria~\eqref{eq:conv_criterion_abs}--\eqref{eq:conv_criterion_rel}}
        \State $\timeIdx \gets \timeIdx + 1$
        \Comment increment time-step index
    \Until{$\timeIdx = \numSnaps$}
    \Comment{$\numSnaps$ is the total number of time steps}
\end{algorithmic}

 \label{alg:schwarz}
\end{algorithm}

\begin{remark}
An important property of the Schwarz alternating method is that it allows, in its most general form, the use of different time integrators as well as different time steps within different subdomains. The assumption that all subdomains are advanced using the same time-step can be relaxed by introducing the idea of a controller time-stepper; please see~\cite{Mota:2022} for details. It is additionally demonstrated in~\cite{Mota:2022} that this application of the Schwarz alternating method is equivalent to performing Schwarz in a space-time domain.
\end{remark}

Note that we have specifically presented the \textit{additive} Schwarz algorithm. This is in contrast to the original \textit{multiplicative} Schwarz algorithm considered in our past works~\cite{Mota:2017, Mota:2022, Mota:2023, Barnett:2022}, which solves each subdomain-local problem in sequence, rather than concurrently, thus incorporating the most up-to-date interface information. To modify Algorithm~\ref{alg:schwarz}, line 8 would simply be altered to interpolate the solution $\solVar^{(\iterIdx)}_1$ onto $\Gamma_2$. Although it enables concurrent subdomain calculations, the additive Schwarz algorithm will in general require more iterations to converge than would the multiplicative algorithm~\cite{Gander:2008}. However, it has been shown that, provided appropriate definition of the operators $\mathcal{B}_{ij}$ (see Table~\ref{table:Schwarz_BCs}), the additive Schwarz iteration~\eqref{eq:generic_schwarz_iter_Omega1} and~\eqref{eq:generic_schwarz_iter_Omega2} will converge to the same solution of~\eqref{eq:generic_pde} as the multiplicative Schwarz iteration. All results shown in this paper utilize the additive Schwarz algorithm to enable concurrent calculations across all subdomains.

\subsection{Nonstandard Schwarz BCs for cell-centered finite volume discretizations} \label{sec:schwarz_CCFV}

For spatial discretizations which define degrees of freedom at the endpoints, vertices, or edges of mesh elements, the application of non-overlapping Schwarz boundary projection operators (as outlined in Table~\ref{table:Schwarz_BCs}) is rather straightforward. Often, this simply involves interpolation on a line segment for a 2D domain, or on a plane for a three-dimensional domain. Finite difference, finite element, and vertex-centered finite volume discretization schemes are readily adapted to non-overlapping Schwarz decompositions as a result. On the other hand, CCFV schemes define degrees of freedom at the centroid of the mesh volumes. For such schemes, boundary conditions can be applied by a ``ghost cell'' formulation, whereby fictitious cells are placed outside the computational domain, adjacent to cells whose faces compose the domain boundaries. These ghost cells need not have any geometry associated with them, though a centroid coordinate may be used for the sake of convenience. The system state in these ghost cells is specified to weakly enforce boundary conditions through the finite volume fluxes. For example, for a cell at an impermeable slip wall, the accompanying ghost cell state would be set to have the same density, pressure, and tangential velocity of the interior cell, but the opposite normal velocity. The fluxes at the boundary face would then be solved exactly the same as for any other cell face, and implicitly enforce zero transport through the wall. An overview of ghost cell boundary conditions for a variety of other physical boundaries can be found in many classical texts on finite volume methods, such as~\cite{LeVeque2002}.

\begin{figure}
    \centering
    {
    \begin{tikzpicture}
        \draw [blue, line width=1.5] (0,0) -- (7.0,0);
        \draw [blue, line width=1.5] (0,-0.15) -- (0,0.15);
        \draw [blue, line width=1.5] (1.75,-0.15) -- (1.75,0.15);
        \draw [blue, line width=1.5] (3.5,-0.15) -- (3.5,0.15);
        \draw [blue, line width=1.5] (5.25,-0.15) -- (5.25,0.15);
        \draw [blue, line width=1.5] (7.0,-0.15) -- (7.0,0.15);
        \draw [blue, line width=1, dashed] (7.0,0) -- (8.75,0);
        \draw [blue, line width=1, dashed] (8.75,-0.15) -- (8.75,0.15);
        \node [blue] at (0.88,0.45) {$~x_{1,1}$};
        \node [blue] at (2.64,0.45) {$~x_{1,2}$};
        \node [blue] at (4.4,0.45) {$~x_{1,3}$};
        \node [blue] at (6.16,0.45) {$~x_{1,4}$};
        \node [blue] at (7.92,0.45) {$~x_{1,g}$};
        \filldraw [blue] (0.88,0) circle (2pt);
        \filldraw [blue] (2.64,0) circle (2pt);
        \filldraw [blue] (4.4,0) circle (2pt);
        \filldraw [blue] (6.16,0) circle (2pt);
        \draw [orange, line width=1.5] (3.5,-0.3) -- (10.5,-0.3);
        \draw [orange, line width=1.5] (3.5,-0.45) -- (3.5,-0.15);
        \draw [orange, line width=1.5] (5.25,-0.45) -- (5.25,-0.15);
        \draw [orange, line width=1.5] (7.0,-0.45) -- (7.0,-0.15);
        \draw [orange, line width=1.5] (8.75,-0.45) -- (8.75,-0.15);
        \draw [orange, line width=1.5] (10.5,-0.45) -- (10.5,-0.15);
        \draw [orange, line width=1, dashed] (1.75,-0.3) -- (3.5,-0.3);
        \draw [orange, line width=1, dashed] (1.75,-0.45) -- (1.75,-0.15);
        \node [orange] at (4.4,-0.75) {$~x_{2,1}$};
        \node [orange] at (6.16,-0.75) {$~x_{2,2}$};
        \node [orange] at (7.92,-0.75) {$~x_{2,3}$};
        \node [orange] at (9.68,-0.75) {$~x_{2,4}$};
        \node [orange] at (2.64,-0.75) {$~x_{2,g}$};
        \filldraw [orange] (4.4,-0.3) circle (2pt);
        \filldraw [orange] (6.16,-0.3) circle (2pt);
        \filldraw [orange] (7.92,-0.3) circle (2pt);
        \filldraw [orange] (9.68,-0.3) circle (2pt);
    \end{tikzpicture}
    \subcaption{Overlapping interface.}
    \label{subfig:fv_ghost_over}
    }
    \vspace{1em}
    {
    \begin{tikzpicture}
        \draw [blue, line width=1.5] (0,0) -- (7.0,0);
        \draw [blue, line width=1.5] (0,-0.15) -- (0,0.15);
        \draw [blue, line width=1.5] (1.75,-0.15) -- (1.75,0.15);
        \draw [blue, line width=1.5] (3.5,-0.15) -- (3.5,0.15);
        \draw [blue, line width=1.5] (5.25,-0.15) -- (5.25,0.15);
        \draw [blue, line width=1.5] (7.0,-0.15) -- (7.0,0.15);
        \draw [blue, line width=1, dashed] (7.0,0) -- (8.75,0);
        \draw [blue, line width=1, dashed] (8.75,-0.15) -- (8.75,0.15);
        \node [blue] at (0.88,0.45) {$~x_{1,1}$};
        \node [blue] at (2.64,0.45) {$~x_{1,2}$};
        \node [blue] at (4.4,0.45) {$~x_{1,3}$};
        \node [blue] at (6.16,0.45) {$~x_{1,4}$};
        \node [blue] at (7.92,0.45) {$~x_{1,g}$};
        \filldraw [blue] (0.88,0) circle (2pt);
        \filldraw [blue] (2.64,0) circle (2pt);
        \filldraw [blue] (4.4,0) circle (2pt);
        \filldraw [blue] (6.16,0) circle (2pt);
        \draw [orange, line width=1.5] (7.0,-0.3) -- (14.0,-0.3);
        \draw [orange, line width=1.5] (7.0,-0.45) -- (7.0,-0.15);
        \draw [orange, line width=1.5] (8.75,-0.45) -- (8.75,-0.15);
        \draw [orange, line width=1.5] (10.5,-0.45) -- (10.5,-0.15);
        \draw [orange, line width=1.5] (12.25,-0.45) -- (12.25,-0.15);
        \draw [orange, line width=1.5] (14.0,-0.45) -- (14.0,-0.15);
        \draw [orange, line width=1, dashed] (5.25,-0.3) -- (7.0,-0.3);
        \draw [orange, line width=1, dashed] (5.25,-0.45) -- (5.25,-0.15);
        \node [orange] at (7.92,-0.75) {$~x_{2,1}$};
        \node [orange] at (9.68,-0.75) {$~x_{2,2}$};
        \node [orange] at (11.44,-0.75) {$~x_{2,3}$};
        \node [orange] at (13.2,-0.75) {$~x_{2,4}$};
        \node [orange] at (6.16,-0.75) {$~x_{2,g}$};
        \filldraw [orange] (7.92,-0.3) circle (2pt);
        \filldraw [orange] (9.68,-0.3) circle (2pt);
        \filldraw [orange] (11.44,-0.3) circle (2pt);
        \filldraw [orange] (13.2,-0.3) circle (2pt);
    \end{tikzpicture}
    \subcaption{``Non-overlapping'' interface}
    \label{subfig:fv_ghost_nonover}
    }
    
    \caption{Example 1D cell-centered finite volume Schwarz interfaces. Physical cells are marked by solid lines, while ghost cells are marked by dashed lines. A spatial offset is included for clarity, but overlapping cells should be considered to exactly align with their neighboring subdomain counterparts.}
    \label{fig:fv_ghost}
\end{figure}
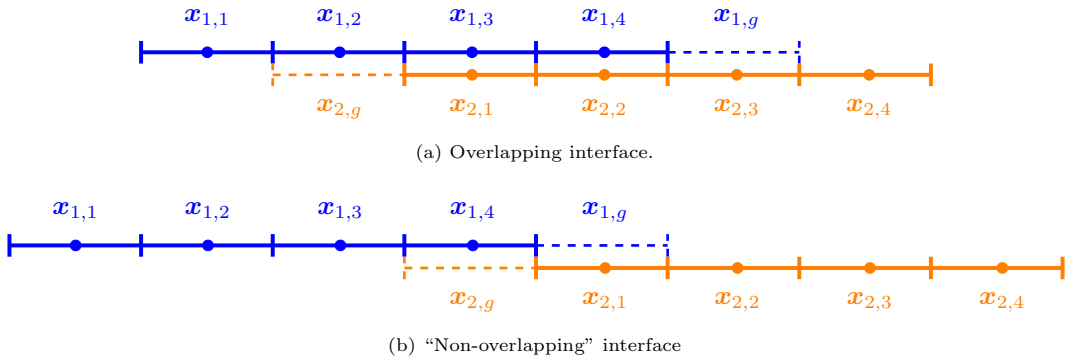

Applying the Schwarz alternating method to a domain decomposition with a CCFV discretization is largely indistinguishable from its application to any other discretization. Figure~\ref{subfig:fv_ghost_over} shows a simplified graphical depiction of a uniform 1D CCFV overlapping decomposition. Solid lines indicate physical domains, while dashed lines indicate ghost cells. The Schwarz projection operator simply interpolates the interior state of one domain to the ghost cells of the neighboring subdomains. Solving the boundary face fluxes in the neighboring subdomain thus accounts for the effects of the solution as it is updated through the Schwarz algorithm.

The use of ghost cells in a finite volume discretization leads to a somewhat nonstandard treatment of non-overlapping Schwarz interfaces, as demonstrated visually in Figure~\ref{subfig:fv_ghost_nonover}. While the physical subdomains $\Omega_1$ and $\Omega_2$ do not overlap and only share a physical cell face at a single point (the midpoint between $x_{1,4}$ and $x_{2,1}$), the ghost cells create an effective overlap region having a width of two cells. This permits the application of a Dirichlet--Dirichlet Schwarz interface condition by assigning the state of the physical cells at the boundary to the corresponding ghost cell of the adjacent subdomain. The result is a more versatile method, as it is possible to obtain a convergent sequence of Schwarz iterations when utilizing a (more flexible) non-overlapping DD with Dirichlet transmission BCs at the Schwarz boundaries (row 1 of Table~\ref{table:Schwarz_BCs}). This greatly simplifies the interface treatment by avoiding Dirichlet--Neumann or Robin--Robin conditions as in traditional non-overlapping Schwarz schemes.

While the scheme described above is not in truth a non-overlapping scheme, as there is an effective overlap region defined by the ghost cells, the resulting Schwarz variant possesses the same flexibility of a ``true'' non-overlapping method by enabling the independent discretization of neighboring subdomains and application to interface problems (see Remark 2). The fact that ghost cells are non-physical constructs hints at interesting possibilities for a simplified Schwarz decomposition for complex finite volume meshes. This treatment is not dissimilar from halo exchange regions utilized in parallel computing, where neighboring parallel partitions must communicate data at interfaces to propagate dynamics throughout the full domain. Unlike parallel partitioning, the Schwarz iterative method provides a simple and rigorous coupling mechanism for subdomains which may be governed by dissimilar physics, time integrators, or mesh topologies. The utility of this ``non-overlapping'' finite volume Dirichlet-Dirichlet interface will be demonstrated empirically in Section~\ref{sec:results}.

\subsection{Schwarz convergence properties for CCFV discretizations} \label{sec:convergence}

The convergence of the Schwarz alternating method has been studied by a number of authors for a variety of different PDEs, e.g.,~\cite{Lions:1988, LionsNonOverlap1988, Zanolli:1987, Gander:2003, Mota:2017, Mota:2022}. Analysis of the method for the specific case of nonlinear conservation laws such as those considered here was examined by Gander and Rohdi~\cite{Gander:2005}, and is not repeated here.

While we are not aware of any publications discussing in general the convergence of the non-overlapping Schwarz alternating method with Dirichlet-Dirichlet transmission BCs for the specific case of CCFV discretizations, it is shown by Dolean \textit{et al.}~\cite{Dolean:2007, Dolean:2009} that the Dirichlet-Dirichlet form of the method can be convergent for certain linear hyperbolic PDEs regardless of the underlying spatial discretization. Attention is restricted to the Cauchy--Riemann~\cite{Dolean:2007} and Maxwell equations~\cite{Dolean:2009}. It is demonstrated that, when discretized in time, an initial boundary value problem (IBVP) having Dirichlet BCs for one of these PDEs is equivalent to an elliptic problem having Robin BCs. Since the Schwarz alternating method is provably convergent for elliptic problems with Robin-Robin transmission BCs~\cite{Zanolli:1987, Gander:2008}, the authors make the argument that the Dirichlet-Dirichlet Schwarz iteration will converge for the hyperbolic IBVPs from which the elliptic PDE was derived. We attempted to repeat the analysis in~\cite{Dolean:2007, Dolean:2009} for the nonlinear hyperbolic conservation laws considered here, but were unable to show the same equivalence between the time-discretized versions of these PDEs and a corresponding elliptic PDE, due largely to the presence of nonlinear terms in the equations of interest. As discussed earlier in Section~\ref{sec:schwarz_CCFV}, we believe that the convergence of the non-overlapping Schwarz iteration with Dirichlet-Dirichlet transmission BCs observed here (see Section~\ref{sec:results}) is due to the small overlap introduced into the discretization via ghost cells used for BC enforcement (Figure~\ref{fig:fv_ghost}). This convergence phenomenon can be explained using an intuitive argument. In a typical node-centered discretization, applying a Dirichlet-Dirichlet transmission BC on a non-overlapping boundary will not allow sufficient propagation of solution information from the interior of one subdomain onto the boundary of its neighboring subdomain. This is not the case, however, in discretizations where the primary degrees of freedom are defined at the cell centers and boundary data are communicated via ghost cells. Here, data from the \textit{interior} of each subdomain (half an element in) are communicated through the Dirichlet BCs imposed on the Schwarz boundaries, and it is for this reason that the Dirichlet-Dirichlet Schwarz iteration is convergent for CCFV discretizations.

It is worth noting that the ghost cells used in CCFV methods introduce only a minimum amount of overlap (Figure~\ref{fig:fv_ghost}). It is well-known~\cite{Lions:1988, Gander:2008, Mota:2022, Tran:2011} that the convergence rate of the overlapping Schwarz alternating method is directly related to the size of the overlap region; faster convergence is expected when a larger overlap region is defined. We thus expect our Dirichlet-Dirichlet Schwarz method to require more iterations to reach convergence when used to couple non-overlapping subdomains discretized using CCFV discretizations (which are really overlapping domains having a minimal amount of overlap). Even with this inefficiency, we are often able to achieve speedups for our PROM-based couplings with respect to a corresponding monolithic FOM, as shown in Section~\ref{sec:results}. This is partially due to the fact that employing a non-overlapping DD eliminates redundant calculations for cells in the overlap region during the Schwarz iteration process. It may be possible to improve convergence of our non-overlapping Schwarz method through the development of optimized transmission BCs~\cite{Gander:2008, Deng:2003, Lui:2001, Gerardo-Giorda:2013, Nataf:1997, Gander:2003, Halpern:2008, Nataf:1994}, which can have a mix of boundary conditions, as discussed earlier in Remark 1. We do not attempt that here, but it may be an interesting extension to explore in future work.

% Projection-based ROM
\section{Projection-based reduced-order models} \label{sec:rom}

The numerical solution of~\eqref{eq:generic_pde}, or its Schwarz formulation in~\eqref{eq:generic_schwarz_iter_Omega1} and~\eqref{eq:generic_schwarz_iter_Omega2}, requires discretization in both time and space. For a general spatial discretization (e.g., finite element, finite volume, etc.), the resulting ordinary differential equation (ODE) can be expressed as

\begin{equation}\label{eq:semidiscrete}
    \ode{\solVec (t)}{t} + \nonlinOpVec \left( \solVec (t) \right) = \forceVec, \quad \solVec(0) = \solVec^0,
\end{equation}
where $\funcMap{\solVec}{[0, T]}{\mathbb{R}^{\numDOFs}}$ is the discrete solution, $\funcMap{\nonlinOpVec}{\mathbb{R}^{\numDOFs}}{\mathbb{R}^{\numDOFs}}$ is the nonlinear vector-valued equivalent of the nonlinear operator $\nonlinOp$, and $\forceVec \inROne{\numDOFs}$ is the source vector. Discretizing~\eqref{eq:semidiscrete} in time with an arbitrary time integrator results in the simplified residual formulation

\begin{equation}\label{eq:residual}
    \resFunc{\solVec^{\timeIdx}} = \mathbf{0},
\end{equation}
where $\funcMap{\resVec}{\mathbb{R}^{\numDOFs}}{\mathbb{R}^{\numDOFs}}$ is the fully-discrete residual, and $\solVec^{\timeIdx} \inROne{\numDOFs}$ is shorthand for the solution at time instance $t_{\timeIdx}$.

For systems which require high spatial resolutions (large $\numDOFs$) to achieve a sufficiently accurate solution, or complex systems for which evaluating the nonlinear terms $\nonlinOpVec(\cdot)$ is expensive, the computational cost of evaluating~\eqref{eq:residual} may be prohibitive. Any method which seeks to limit this expense by drastically reducing the number of degrees of freedom $\numDOFs$ associated with solving~\eqref{eq:residual}, with minimal loss of accuracy, can be generally referred to as a model order reduction method. In this section, we describe the particular approach of PROMs and their application to Schwarz-based domain decomposition.

\subsection{Overview of PROMs}

To begin, PROMs propose an approximation of the solution $\solVec$ which is represented by a relatively small number $\numModes \ll \numDOFs$ of degrees of freedom. While nonlinear approximations have been proposed~\cite{Lee:2020}, we restrict our discussion to linear approximations of the form

\begin{equation}\label{eq:approx_sol}
    \solVec (t) \approx \solVecRom (t) \defEq \solVecCent + \trialBasis \solVecCoef (t),
\end{equation}
where $\funcMap{\solVecRom}{[0, T]}{\numDOFs}$ is the approximate state, and $\solVecCent \inROne{\numDOFs}$ is a centering (or reference) state. The trial basis $\trialBasis \defEq [\trialVec_1, \ \hdots, \ \trialVec_\numModes] \inRTwo{\numDOFs}{\numModes}$ is composed of $\numModes$ trial basis vectors $\trialVec_i$, while $\funcMap{\solVecCoef}{[0, T]}{\numModes}$ are the generalized coordinates (or modal coefficients) associated with the trial basis.

Inserting the approximate solution in~\eqref{eq:approx_sol} into the governing ODE~\eqref{eq:semidiscrete} results in the formulation

\begin{equation}\label{eq:romN}
    \trialBasis \ode{\solVecCoef}{t} + \nonlinOpVec(\solVecCent + \trialBasis \solVecCoef) = \forceVec,
\end{equation}
which is still an $\numDOFs$-dimensional ODE. Projection by an appropriate test basis $\testBasis \inRTwo{\numDOFs}{\numModes}$ and rearranging terms arrives at the desired reduced-order model

\begin{equation}\label{eq:rom_ode}
    \ode{\solVecCoef}{t} + \projMat \nonlinOpVec(\solVecCent + \trialBasis \solVecCoef) = \projMat \forceVec,
\end{equation}
where we define the term $\projMat \defEq [\testBasis^\top \trialBasis]^{-1} \testBasis^\top \inRTwo{\numModes}{\numDOFs}$ for notational simplicity. Ignoring for the moment the fact that evaluating the non-linear function still scales with $\numDOFs$,~\eqref{eq:rom_ode} is an $\numModes$-dimensional ODE, and its integration in time should incur a much lower computational cost than that of~\eqref{eq:semidiscrete}. Naturally, the quality of the state approximation in~\eqref{eq:approx_sol} and the robustness of the PROM ODE in~\eqref{eq:rom_ode} are highly dependent on the choice of trial basis $\trialBasis$ and test basis $\testBasis$. We address each in turn.

\subsubsection{Proper orthogonal decomposition (POD)} \label{sec:POD}

While there exist a number of analytical bases which may be used for state approximations (e.g., Fourier series), data-driven alternatives offer a means of tailoring the choice of trial basis $\trialBasis$ to a given solution manifold, theoretically achieving much more accurate representations with far fewer basis vectors. To date, POD~\cite{Holmes:1996, Sirovich:1987} remains the most popular method for deriving such an approximation from data, and has similarly seen wide application under the titles principal component analysis in statistics and the Karhunen-Lo\`{e}ve expansion in mathematics.

POD begins by collecting $\numSnaps$ data ``snapshots'' from a number of high-fidelity simulations (i.e., $\numDOFs$-dimensional, without dimension reduction) into a snapshot matrix $\solMat \defEq [\solVec_1, \ \hdots \ , \ \solVec_{\numSnaps}] \inRTwo{\numDOFs}{\numSnaps}$. These snapshots may be taken from high-fidelity simulations under varying parameterizations and at different time instances. The rank $\numModes$ POD basis is defined as the solution to the least-squares problem,

\begin{equation}\label{eq:pod_min}
    \trialBasis = \underset{\dummyBasis \inRTwo{\numDOFs}{\numModes}}{\text{arg\,min}} \left\Vert \solMat - \dummyBasis \dummyBasis^\top \solMat \right\Vert_2.
\end{equation}
That is, POD seeks the matrix $\trialBasis$ which minimizes the $\ell^2$ error in the orthogonal projection of the data matrix. Equation~\eqref{eq:pod_min} has a simple solution from the singular value decomposition, whereby the data matrix is decomposed as $\solMat = \leftSVecMat \EValMat \rightSVecMat^\top$, with $\leftSVecMat \inRTwo{\numDOFs}{\numSnaps}$, $\EValMat \inRTwo{\numSnaps}{\numSnaps}$, and $\rightSVecMat \inRTwo{\numSnaps}{\numSnaps}$ being the left singular vectors, singular value diagonal matrix, and right singular vectors respectively (under the assumption $\numDOFs \ge \numSnaps$). The orthonormal trial basis $\trialBasis$ is simply computed as the leading $\numModes$ leading left singular vectors of $\leftSVecMat$.

Increasing the number of basis vectors $\numModes$ retained in the trial basis $\trialBasis$ will monotonically decrease the data projection error associated with~\eqref{eq:pod_min}. However, there are no such guarantees that increasing $\numModes$ will improve out-of-sample projection error (i.e., for snapshots taken from time instances or parameter values not included in $\solMat$), or decrease ``online'' error in evaluating the PROM ODE~\eqref{eq:rom_ode} due to the nonlinear treatment of errors. In fact, increased $\numModes$ for fields characterized by strong gradients often exacerbates oscillatory ``ringing'' phenomena inherent in linear representations, which may destabilize a PROM simulation. Methods exist for limiting such behavior by artificial dissipation or filtering~\cite{Ahmed2021,Xie2018}, and others exist to adapt the trial basis during the PROM runtime~\cite{Peherstorfer:2020} to better represent the instantaneous solution. Here, we restrict analysis to static trial bases and do not explicitly address spurious oscillatory behavior. When applying PROMs to nonlinear systems of sufficient complexity, the basis dimension is usually not the largest driver of computational cost (as will be discussed in Section~\ref{sec:hr}), and the number of modes may generally be set large enough to accurately represent the solution.

\subsubsection{Least-squares Petrov--Galerkin (LSPG) projection} 

In selecting the projecting test basis $\testBasis$, by far the simplest choice is Galerkin projection, which simply assigns $\testBasis = \trialBasis$. Any other choice of $\testBasis$ is referred to as as Petrov--Galerkin projection. Galerkin projection has seen widespread success for parabolic and elliptic PDEs~\cite{Quarteroni2011}, but it often leads to unstable solutions for advection-dominated hyperbolic systems~\cite{Grimberg2020}. Although it has been shown that recasting the governing equations with a particular alternate set of variables and applying Galerkin projection may resolve these issues~\cite{Parish2023}, the dominant alternative to Galerkin projection is the least-squares Petrov--Galerkin (LSPG) projection developed by Carlberg \textit{et al.}~\cite{Carlberg:2011}. 

Unlike Galerkin projection, which can be considered the minimization of the semidiscrete PROM ODE residual in~\eqref{eq:rom_ode}, LSPG is constructed under the goal of minimizing the fully-discrete residual in~\eqref{eq:residual}. At a given time step $\timeIdx$, this can be formalized as

\begin{equation}\label{eq:lspg_res}
    \solVecCoef^{\iterIdx} = \underset{\dummyVec \inROne{\numModes}}{\text{arg\,min}} \left\Vert \resOp \resFunc{\solVecCent + \trialBasis \dummyVec} \right\Vert_2^2,
\end{equation}
where the standard LSPG method takes $\resOp = \mathbf{I}_{\numDOFs}$. The solution of~\eqref{eq:lspg_res} by the Gauss--Newton algorithm arrives at the iterative solution

\begin{equation}\label{eq:lspg_normal}
    \left[ \jacobMat^{\iterIdx} \trialBasis \right]^\top \jacobMat^{\iterIdx} \trialBasis \Delta \solVecCoef^{\iterIdx} = - \left[ \jacobMat^{\iterIdx} \trialBasis \right]^\top \resFunc{\solVecCent + \trialBasis \solVecCoef^{\iterIdx}},
\end{equation}
where $\jacobMat^{\iterIdx} \defEq \partial \resFunc{\solVecRom^{\iterIdx}} / \partial \solVec $ is the residual Jacobian at iteration $\iterIdx$ of a particular time step, and $\Delta \solVecCoef^{\iterIdx} \defEq \solVecCoef^{\iterIdx+1} - \solVecCoef^{\iterIdx}$ is the solution update. Equation~\eqref{eq:lspg_normal} thus has the form of a Petrov--Galerkin projection with test basis $\testBasis = \jacobMat^{\iterIdx} \trialBasis$.

While LSPG presents additional computational challenges relative to Galerkin projection (namely a time-varying test basis), it has exhibited superior accuracy and robustness when applied to challenging advection-dominated fluid flows~\cite{Carlberg:2017,Grimberg2020}. All results presented in this paper will utilize LSPG in their PROM solutions.
 
\subsubsection{Hyper-reduction} \label{sec:hr}

For general nonlinear systems, the solution of~\eqref{eq:lspg_normal} still scales with the full-order dimension $\numDOFs$ due to the nonlinear term $\nonlinOpVec(\cdot)$. Combined with forming the test basis $\testBasis$ and computing projections, the cost of the PROM is often \textit{greater} than that of the high-fidelity model. To overcome this issue, ``hyper-reduction'' methods seek to eliminate this scaling with $\numDOFs$. We refer to any PROM which implements such a hyper-reduction approach as an HPROM. 

\begin{figure}
    \centering
    \includegraphics[width=0.6\linewidth]{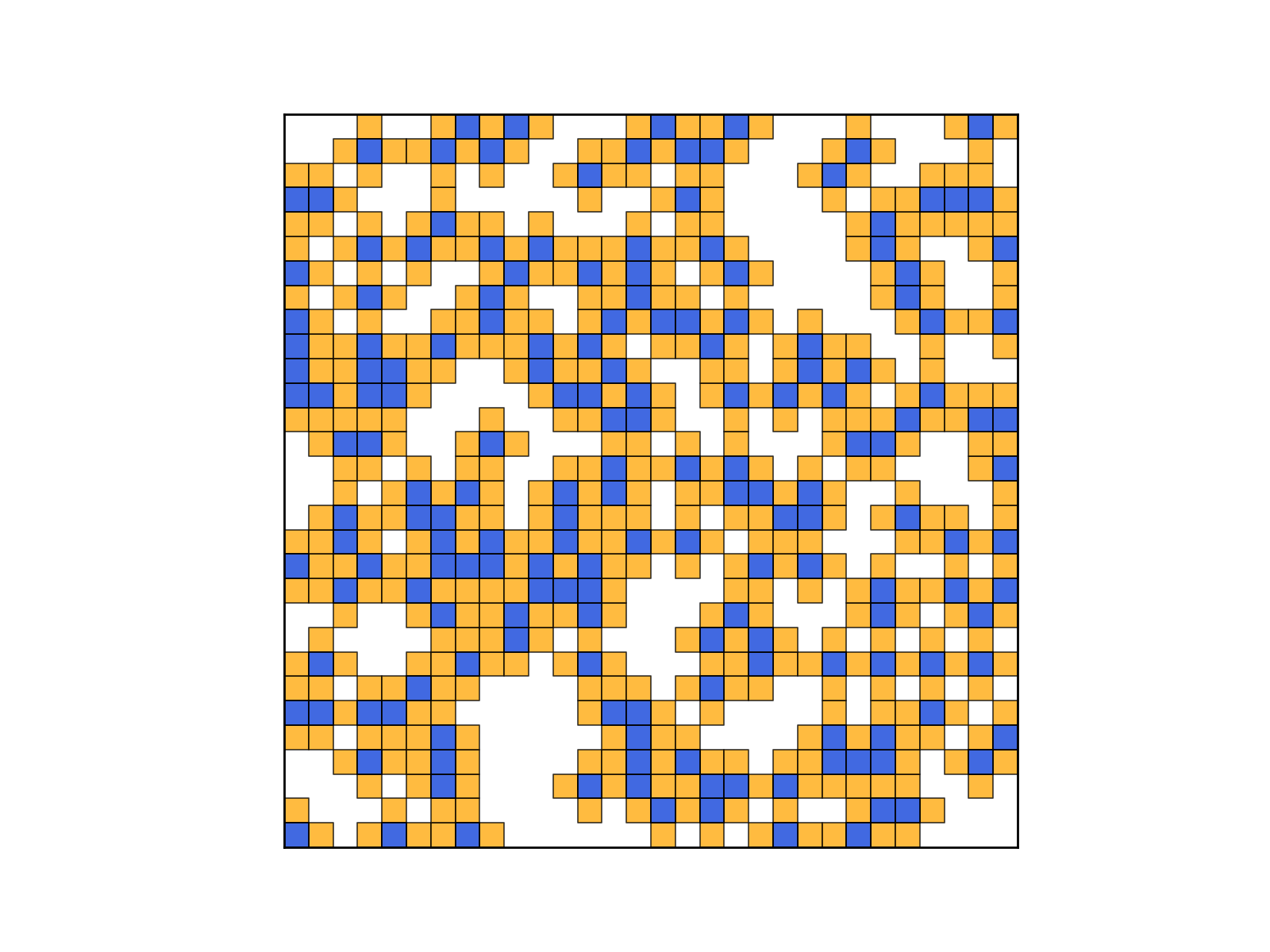}
    \caption{Example sample mesh on a 2D square domain. Blue elements indicate elements where the residual is evaluated, while orange elements indicate supplementary elements required for computing gradients, fluxes, etc.}
    \label{fig:sample_mesh_example}
\end{figure}

Most nonlinear hyper-reduction methods are rooted in the concept of a ``sample mesh,'' whereby the cost of evaluating the nonlinear function $\nonlinOpVec(\cdot)$ is reduced by only evaluating it at a small subset $\numSamps$ of mesh locations. An example of such a sample mesh is displayed in Figure~\ref{fig:sample_mesh_example}. Note the distinction of two different colored elements: some elements are strictly required for evaluating the nonlinear terms at the given subset of mesh locations (shown in blue in Figure~\ref{fig:sample_mesh_example}), while a number of neighboring cells may be required to evaluate these terms as a result of the spatial discretization of the governing PDE (shown in orange in Figure~\ref{fig:sample_mesh_example}). Together, these mesh elements represent the sample mesh, and decreasing the size of this sample mesh not only decreases the computational cost of evaluating the PROM ODE, but also the memory requirements for these calculations. In the remainder of this subsection, we will utilize the ``sampling operator'' $\sampMat = [\standVec_{i_1}, \ \hdots \ , \ \standVec_{i_{\numSamps}}]^\top \inRTwo{\numSamps}{\numDOFs}$, where $e_i$ is the $i$th standard basis vector. Thus the operation $\sampMat \nonlinOpVec(\cdot)$ has the effect of selecting the $\mathcal{S} = \{ i_1, \ \hdots \ , \ i_{\numSamps} \}$ indices of the nonlinear function. In practical calculations, this implies calculating the function only for the elements $\mathcal{S}$ using the degrees of freedom associated with the sample mesh. Note that for multivariate systems, the sample indices do not refer to single mesh elements, but rather individual degrees of freedom (of which a single mesh element may have several).

Different hyper-reduction methods attempt to mitigate the degraded HPROM accuracy and stability associated with approximating the effect of the sample mesh dynamics on the full solution. All results in this paper will utilize the collocation method. Collocation, when applied to LSPG PROMs, takes the extremely simple approach of setting $\resOp = \sampMat$ in~\eqref{eq:lspg_res}, with the resulting normal equations being

\begin{equation}\label{eq:lspg_colloc}
    \left[ \sampMat \jacobMat^{\iterIdx} \trialBasis \right]^\top \sampMat \jacobMat^{\iterIdx} \trialBasis \Delta \solVecCoef^{\iterIdx} = - \left[ \sampMat \jacobMat^{\iterIdx} \trialBasis \right]^\top \sampMat \resFunc{\solVecCent + \trialBasis \solVecCoef^{\iterIdx}}.
\end{equation}
Evaluating~\eqref{eq:lspg_colloc} thus only requires calculating $\numSamps$ rows of the nonlinear residual $\resFunc{\cdot}$ and its Jacobian $\jacobMat^{\iterIdx}$, eliminating any scaling with the full dimension $\numDOFs$. Because collocation is relatively simplistic, it may suffer from accuracy issues by entirely neglecting the effects of unsampled degrees of freedom. This paper does not present any comparisons between these methods or other hyper-reduction methods; the interested reader is directed to other works for such analyses~\cite{Farhat:2015}.

Choosing appropriate sample mesh points comprising a limited budget of  $\numSamps$ total sample mesh points is often of critical importance in ensuring a robust and accurate HPROM solution. For certain applications, the use of simple randomized sampling is sufficient. For other applications in which spatially-localized dynamics may be highly transient, experience strong gradients, or broadly affect other areas of the domain, careful selection of sample points often leads to improved solutions. As finding an optimal sample set is computationally intractable, a wide variety of works propose greedy sampling methods~\cite{Carlberg2013, Peherstorfer2020, Zimmermann2024}. Others offer means of dynamically updating the sample mesh to better capture transient dynamics~\cite{Peherstorfer2015}. For the problems investigated in this paper, for which waves traverse nearly the entirety of the domain, we find that a random sampling approach performs similarly to several greedy approaches, and opt for random sampling for the sake of simplicity. However, the 2D Burgers' equation case benefits from seeding the sample set with those elements selected by the QR decomposition of the trial basis, known as the QDEIM algorithm (related to the discrete empirical interpolation method, or DEIM)~\cite{Drmac2016}, followed by random sampling.

\subsection{Boundary sampling for HPROM coupling}\label{subsec:bound_samp}

As noted in the literature~\cite{Carlberg2013}, deliberately sampling mesh elements at boundaries promotes the incorporation of boundary condition effects on HPROM dynamics. While one may reasonably expect that certain boundary conditions permit zero sampling (e.g., supersonic outflows) or low sampling rates (e.g., far field boundaries), Schwarz interfaces which transmit complex interior dynamics often require careful sampling treatments. In the case of extremely sparse or non-existent sampling at Schwarz interfaces, empirical evidence (not presented here) suggests that the inability to accurately transmit interface information invariably leads to an unstable solution. This follows logically from the spatially-local nature of conservation laws, resulting from an inability to propagate waves across boundaries. 

\begin{figure}
    \begin{minipage}{0.32\linewidth}
        \includegraphics[width=0.99\linewidth]{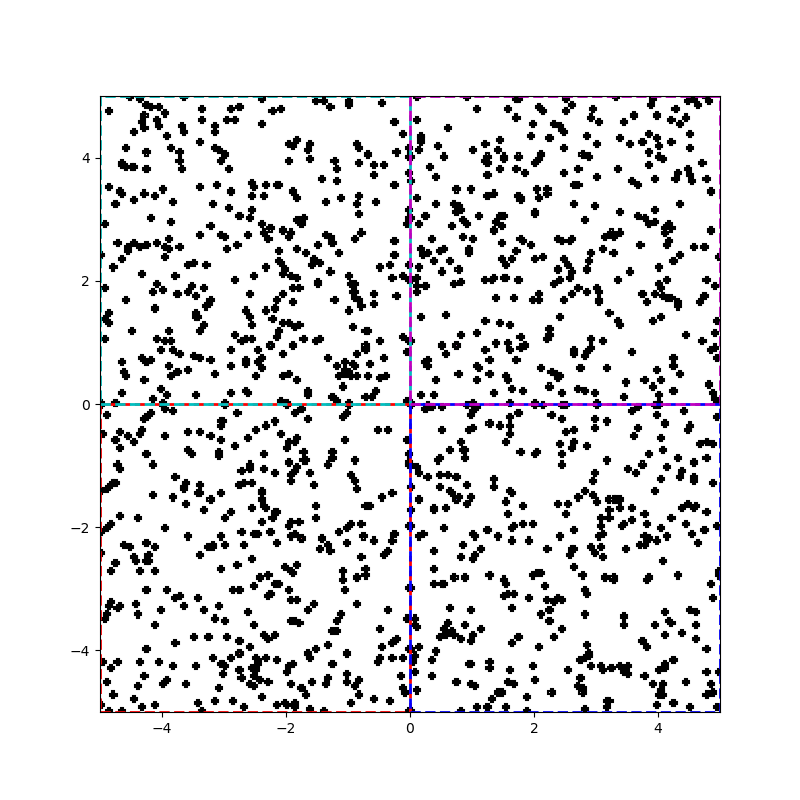}
        \subcaption{$\numBound = 30$}
    \end{minipage}
    \begin{minipage}{0.32\linewidth}
        \includegraphics[width=0.99\linewidth]{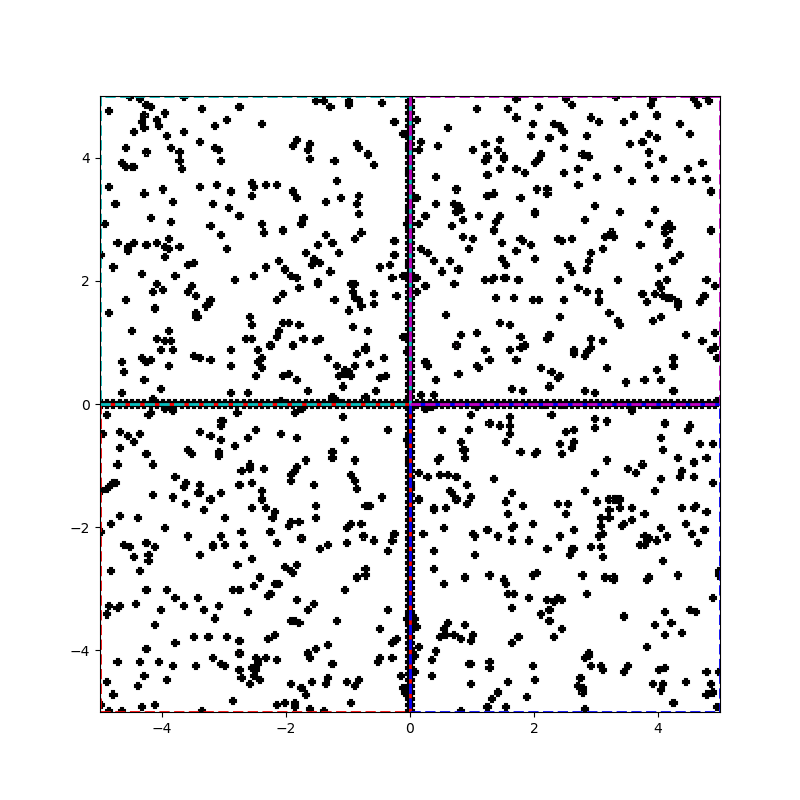}
        \subcaption{$\numBound = 3$}
    \end{minipage}
    \begin{minipage}{0.32\linewidth}
        \includegraphics[width=0.99\linewidth]{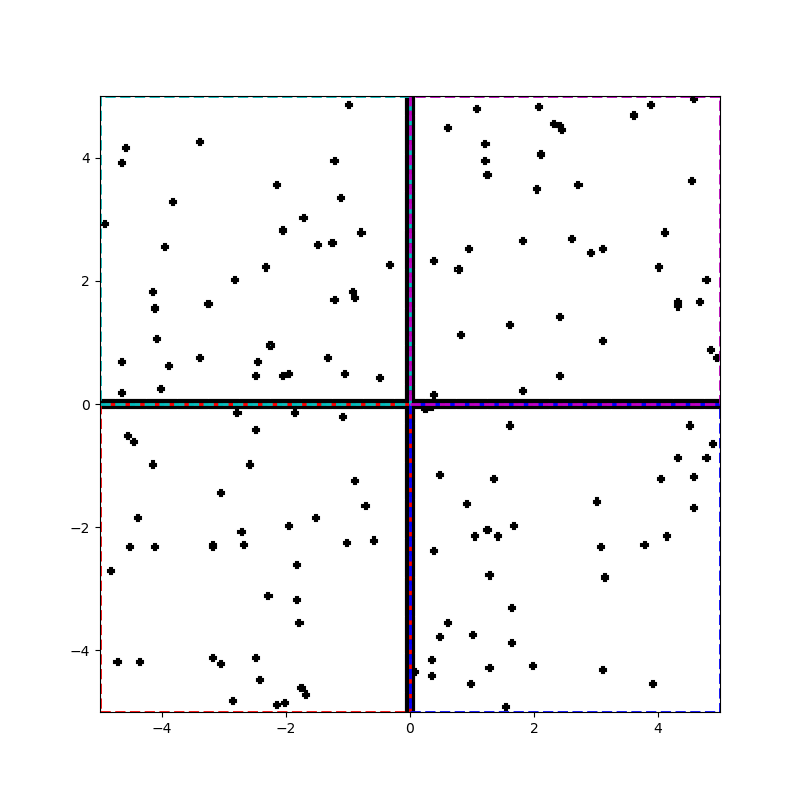}
        \subcaption{$\numBound = 1$}
    \end{minipage}
    \caption{Sample meshes (marked by black dots) for $300 \times 300$ mesh under a $2 \times 2$ domain decomposition (marked by dashed lines), with fixed $\numSamps = 1.5\% \times \numDOFs$ over various interface sampling intervals $\numBound$. Note that a decreased interface sampling interval draws points away from the subdomain interiors.}
    \label{fig:bound_samp_example}
\end{figure}

In this paper, we address this issue with the fairly naive solution of sampling all Schwarz interfaces at a fixed interval $\numBound$. For example, if $\numBound = 5$, then every fifth cell at the boundary is sampled (separate from any additional random or greedy sampling). In the absence of \textit{a priori} system knowledge which might inform more intelligent sampling (e.g., knowing where a shock forms), this is a simplistic yet effective sampling methodology. However, as will be shown in Section~\ref{sec:results}, the interface sampling interval $\numBound$ must be balanced under a fixed sampling budget $\numSamps$ to ensure both robust interface transmission and interior solution accuracy. That is, decreasing the sampling interval necessarily draws sample points from the interior of each subdomain, often decreasing the accuracy of the unsteady solution. This behavior is illustrated in Figure~\ref{fig:bound_samp_example}.

\subsection{PROM and HPROM coupling via the alternating Schwarz method}

The versatility of the Schwarz alternating method lies in its minimal intrusiveness, relying only on a means of projecting solution information from one subdomain onto any interfaces with neighboring subdomains. This flexibility extends to coupling PROMs and HPROMs with very little additional machinery. In the case of PROMs for nonlinear systems, the approximate full order solution $\solVecRom$ must already be accessible at every Schwarz sub-iteration in order to compute the nonlinear terms $\nonlinOpVec(\cdot)$. The application of the projection operators $P_{ij}$ (detailed previously in Table~\ref{table:Schwarz_BCs}) is thus trivial, as they are unchanged from an equivalent FOM-FOM coupling scenario.

Coupling HPROMs requires slightly more care. Recall the concept of the sample mesh described in Section~\ref{sec:hr}: any sampled degrees of freedom require access to adjacent mesh elements to compute the fully-discrete residual. Sample mesh elements located at a Schwarz interface thus require access to mesh elements in neighboring subdomains which may not ``exist,'' insofar as they are not included in the neighboring subdomain's sample mesh set. Especially for large-scale applications with memory restrictions, any elements not explicitly included in the sample mesh will only be accessible in post-processing. Accommodating these additional neighboring degrees of freedom simply requires appending their indices to the sample set, updating the solution at those elements throughout the Schwarz algorithm, and communicating it to neighboring subdomains appropriately. These operations, and the modifications to the projection operators $P_{ij}$ acting on the sampling operator $\sampMat$, largely only require some additional accounting steps. Otherwise, the Schwarz algorithm is unchanged.

As the Schwarz algorithm is fairly agnostic to the choice of model in each subdomain, arbitrary couplings between FOM, PROM, and HPROM subdomains are possible. However, the goal of reduced-order modeling is minimizing the cost of evaluating a given system. Any decomposed system which includes \textit{any} FOM subdomains (or PROM subdomains for nonlinear systems) will therefore be extremely limited in its ability to achieve any computational speedups. This has been previously noted in~\cite{Barnett:2022}, where FOM-PROM or FOM-HPROM coupling fails to achieve significant runtime reductions. In this work, we will restrict analysis to PROM-PROM coupling (as an intermediate step) and HPROM-HPROM coupling.

% Numerical results
\section{Numerical results} \label{sec:results}

This section presents numerical evidence for the effectiveness of the Schwarz alternating method in coupling PROMs. This is demonstrated for three 2D nonlinear hyperbolic fluid flow systems that exhibit phenomena challenging for both efficient model order reduction and robust coupling. The relevant PDEs are given below in roughly increasing order of modeling difficulty.
\begin{enumerate}
    \item The 2D shallow water equations (SWE):
    \begin{align}\begin{split}
        \pde{h}{t} + \pde{hu}{x} + \pde{hv}{y} &= 0 \\
        \pde{(hu)}{t} + \pde{}{x} \left( hu^2 + \frac{1}{2} gh^2 \right) + \pde{(huv)}{y} &= -\mu v \\
        \pde{(hv)}{t} + \pde{(hvu)}{x} + \pde{}{y} \left( hv^2 + \frac{1}{2} gh^2 \right) &= \mu u
    \end{split}\end{align}
    \item The 2D viscous Burgers' equation:
    \begin{align}\begin{split}
        \pde{u}{t} + \frac{1}{2} \left( \pde{u^2}{x} + \pde{uv}{y} \right) &= D \left( \pdeTwo{u}{x} + \pdeTwo{u}{y} \right) \\
        \pde{v}{t} + \frac{1}{2} \left( \pde{uv}{x} + \pde{v^2}{y} \right) &= D \left( \pdeTwo{v}{x} + \pdeTwo{v}{y} \right)
    \end{split}\end{align}
    \item The 2D Euler equations:
    \begin{align}\begin{split}
        \pde{\rho}{t} + \pde{\rho u}{x} + \pde{\rho v}{y} &= 0 \\
        \pde{(\rho u)}{t} + \pde{}{x} \left( \rho u^2 + p \right) + \pde{(\rho u v)}{y} &= 0 \\
        \pde{(\rho v)}{t} + \pde{(\rho v u)}{x} + \pde{}{y} \left( \rho v^2 + p \right) &= 0 \\
        \pde{(\rho E)}{t} + \pde{}{x} \left( (E + p) u \right) + \pde{}{y} \left( (E + p) v \right) &= 0 \\
    \end{split}\end{align}
\end{enumerate}
For all systems, $u$ and $v$ are the $x$- and $y$-velocity, respectively. For the shallow water equations, $h$ is the water column height, and $g = 9.8$ is the acceleration due to gravity. For the Euler equations, $\rho$ is the fluid density, $p$ is the pressure, $E$ is the energy, and the system is closed by the state equation,
\begin{equation}
    p = (\gamma - 1)\left( \rho E - \frac{1}{2} \rho \left(u^2 + v^2 \right) \right),
\end{equation}
where $\gamma = 1.4$ is the ratio of specific heats.

The space-time domains on which these systems are defined are given in Table~\ref{tab:sim_settings}. All systems are discretized by a cell-centered finite volume scheme on a 300$\times$300-cell uniform Cartesian mesh. The advective fluxes are computed by a first-order Roe scheme, and diffusive fluxes are similarly computed by a first-order scheme. Higher-order schemes are perfectly compatible with CCFV Schwarz coupling, simply by projecting the interior solution onto any additional ghost cells of neighboring subdomains. All cases are integrated in time by the first-order implicit Euler scheme with a fixed time step as defined in Table~\ref{tab:sim_settings}. The final solution time for the Burgers' and Euler equations are set such that the cases experiencing the fastest shock speeds ($D = 1 \times 10^{-4}$ and $p_4 = 1.5$, respectively) roughly traverse the entire domain. The SWE final time is chosen somewhat arbitrarily such that the propagating waves traverse the domain multiple times. The time step sizes are chosen to ensure time integrator stability. While Schwarz-based coupling may utilize explicit time integrators, as in~\cite{Mota:2022}, when the time step size in each subdomain is identical then the Schwarz coupling algorithm is unnecessary and is guaranteed to converge in two iterations by definition. The discretization and accompanying solver are implemented using \texttt{pressio} and \texttt{pressio-demoapps}~\cite{Rizzi2020}\footnote{Available at: \url{https://github.com/Pressio}.}.

The initial and boundary conditions for each system are specified as follows. The SWE are initialized with zero velocity, and water height given by a Gaussian pulse centered on $(x, \ y) = (1, \ 1)$ as defined by $h^0 = 1 + 0.125 \exp (-(x - 1)^2 - (y - 1)^2 )$. Impermeable slip walls are enforced at every physical boundary. The Burgers' equations initial conditions are similarly defined as a Gaussian pulse centered on $(x, \ y) = (-0.5, \ -0.4)$ as $u^0 = v^0 = 0.5 \exp((-(x + 0.5)^2 - (y + 0.4)^2)/0.075)$. Homogeneous Dirichlet conditions are enforced at the left and bottom physical boundaries, while homogeneous Neumann conditions are enforced at the top and right boundaries. The Euler equations initial condition is given as a classical Riemann problem, dividing the domain into four quadrants defined by a vertical dividing line at $x = 0.8$ and a horizontal dividing line at $y = 0.8$. The density, velocity, and pressure in the top-right quadrant are specified along with the pressure in the lower-left quadrant. The unknown densities, velocities, and pressures are computed according to the compatibility relations described by Configuration 3 in~\cite{Schulz-Rinne1993} to ensure the initial conditions generate four shocks traversing either in the negative $x$- or $y$-direction. Here, the velocity in the upper-right quadrant is set to zero, the density to 1.5, and the pressure is used to parameterize the problem. The pressure in the lower-left quadrant is fixed at 0.029. All PROMs for the Euler equations are initialized from the FOM solution at $t = 0.05$, as linear subspaces are generally unable to accurately represent several jump discontinuities without producing an unstable solution. This ``warm-start'' permits a slightly smoother (though still extremely sharp) solution. Homogeneous Neumann conditions are enforced at all physical boundaries.

\begin{table}
    \caption{Numerical simulation settings for each system. $T$ and $\Delta t$ are, respectively, the total simulation time and the time step size.}
    \label{tab:sim_settings}
    \centering
    \begin{tabular}{lrrr}
    \toprule
    Setting & 2D SWE & 2D Burgers' & 2D Euler \\
    \midrule
    $x$ & $[-5, 5]$ & $[-1, 1]$ & $[0, 1]$ \\ 
    $y$ & $[-5, 5]$ & $[-1, 1]$ & $[0, 1]$ \\
    $T$ & 10.0 & 7.5 & 0.9 \\
    $\Delta t$ & 0.01 & 0.05 & 0.005 \\
    \bottomrule
    \end{tabular}
\end{table}

\begin{figure}
    {
        \begin{minipage}{0.32\linewidth}
            \includegraphics[width=0.99\linewidth]{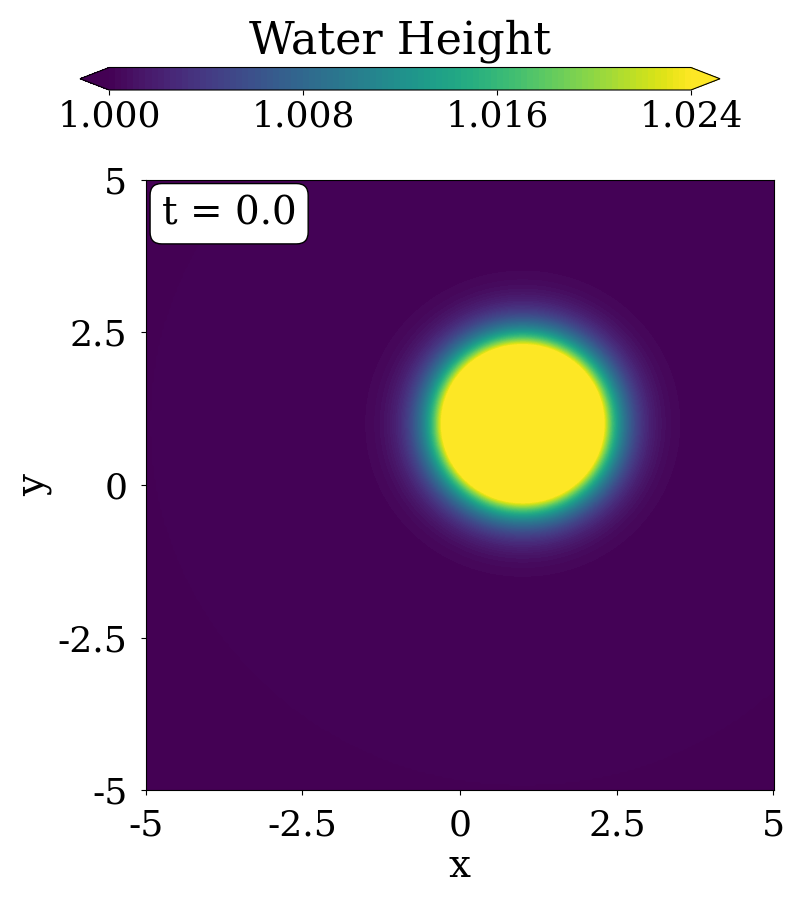}
        \end{minipage}
        \begin{minipage}{0.32\linewidth}
            \includegraphics[width=0.99\linewidth]{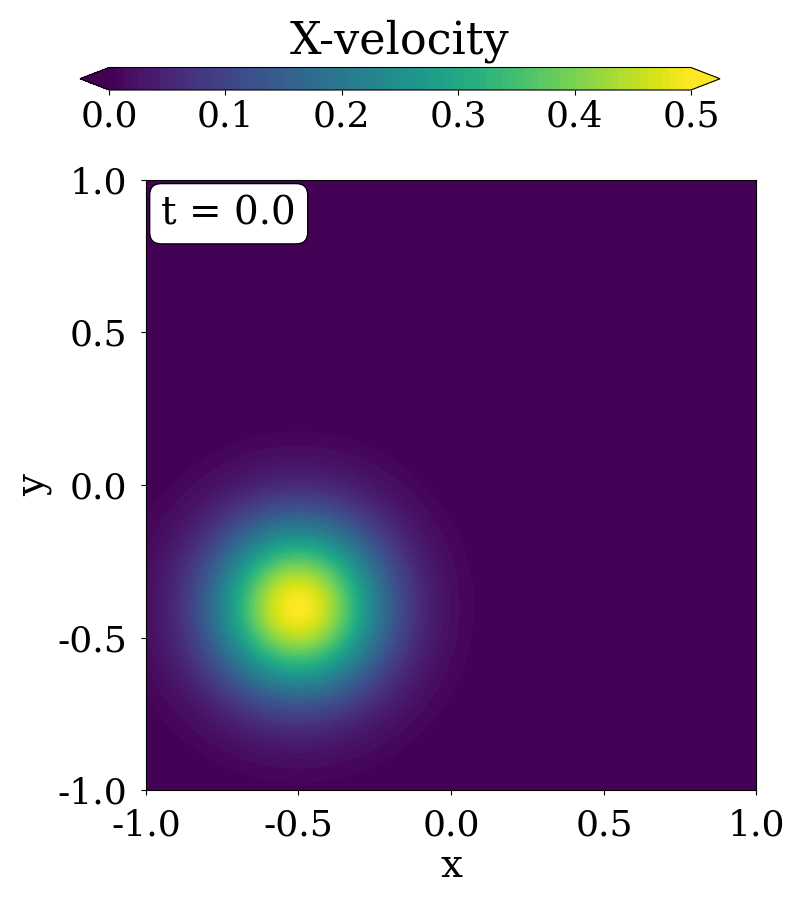}
        \end{minipage}
        \begin{minipage}{0.32\linewidth}
            \includegraphics[width=0.99\linewidth]{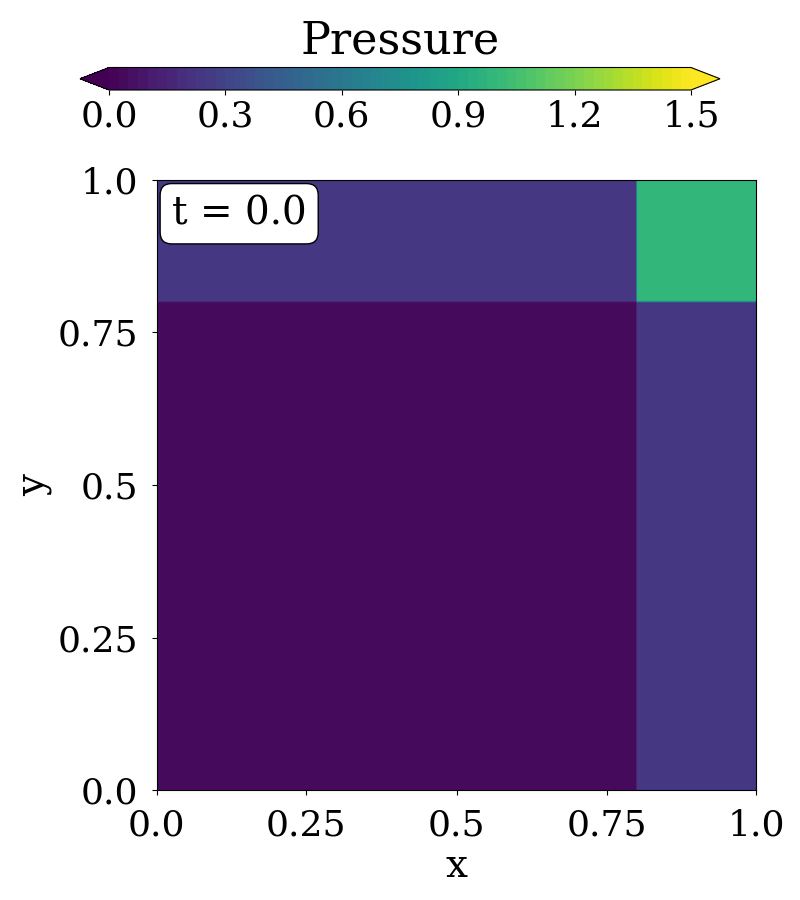}
        \end{minipage}
    }

    {
        \begin{minipage}{0.32\linewidth}
            \includegraphics[width=0.99\linewidth]{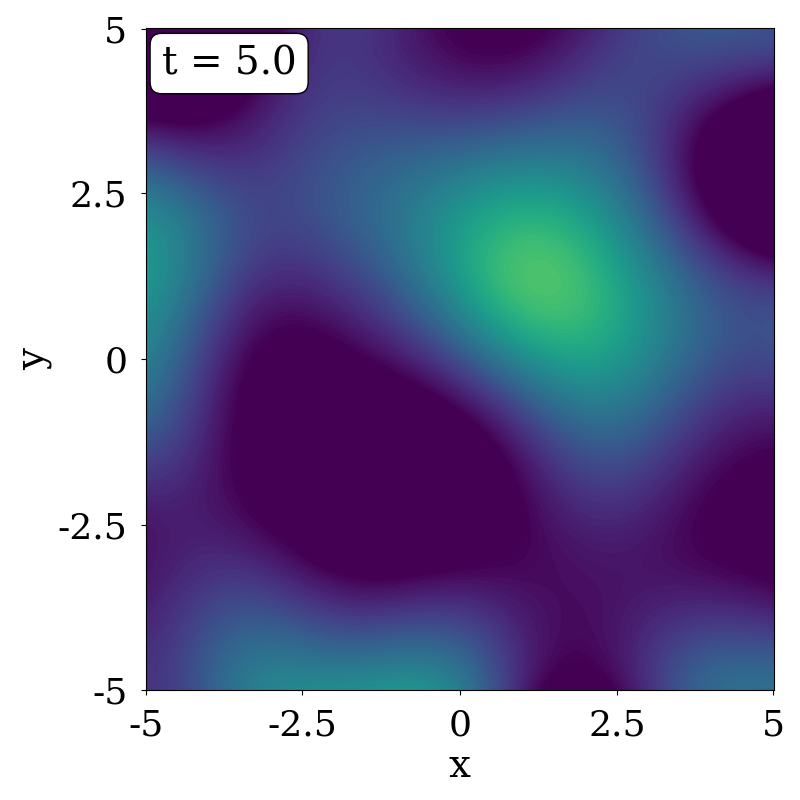}
        \end{minipage}
        \begin{minipage}{0.32\linewidth}
            \includegraphics[width=0.99\linewidth]{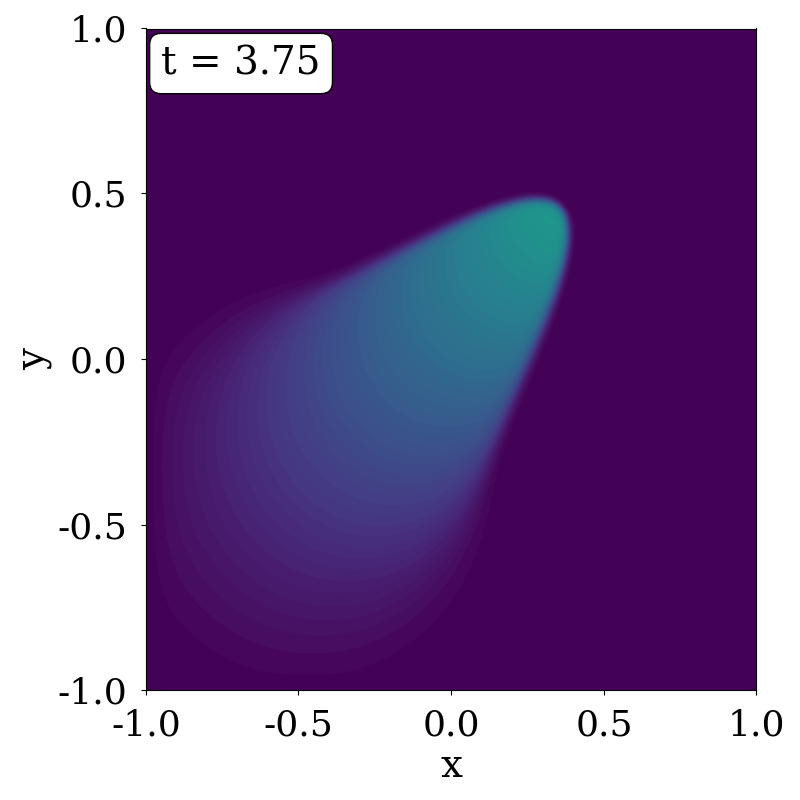}
        \end{minipage}
        \begin{minipage}{0.32\linewidth}
            \includegraphics[width=0.99\linewidth]{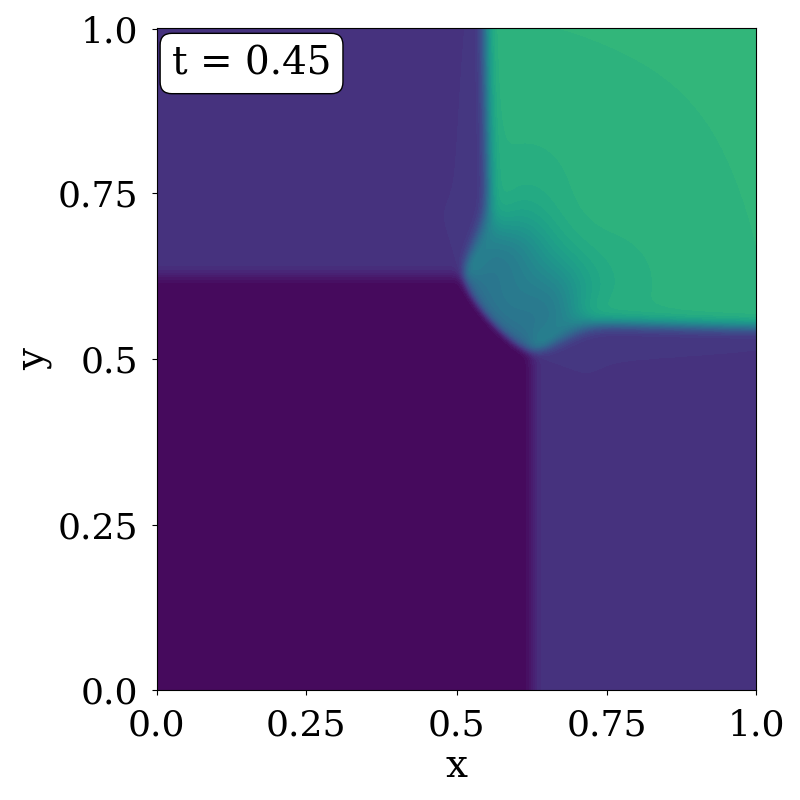}
        \end{minipage}
    }

    {
        \begin{minipage}{0.32\linewidth}
            \includegraphics[width=0.99\linewidth]{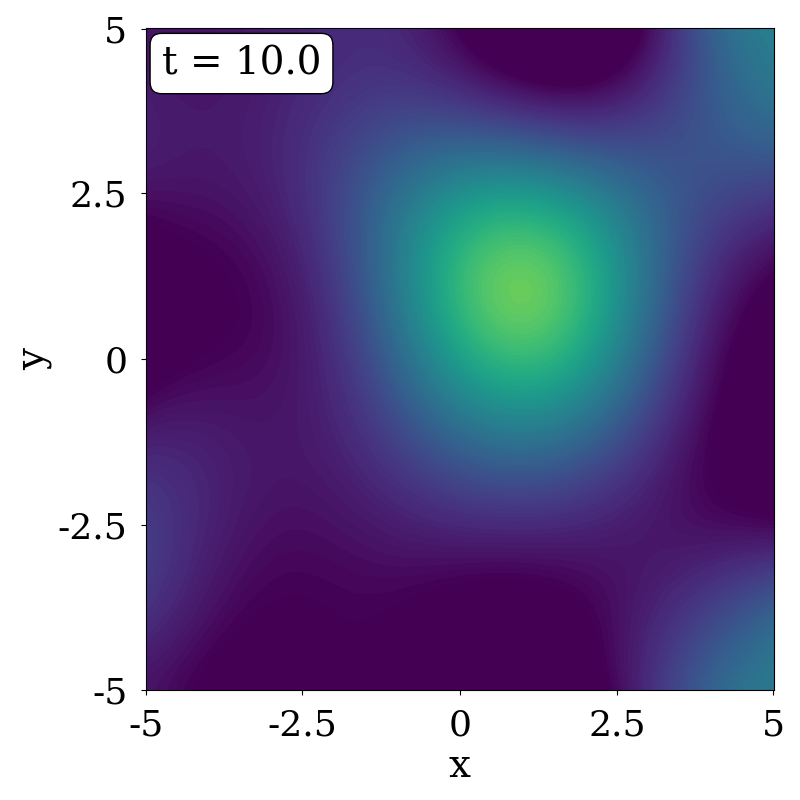}
            \subcaption{SWE, $\mu = -2.0$}
        \end{minipage}
        \begin{minipage}{0.32\linewidth}
            \includegraphics[width=0.99\linewidth]{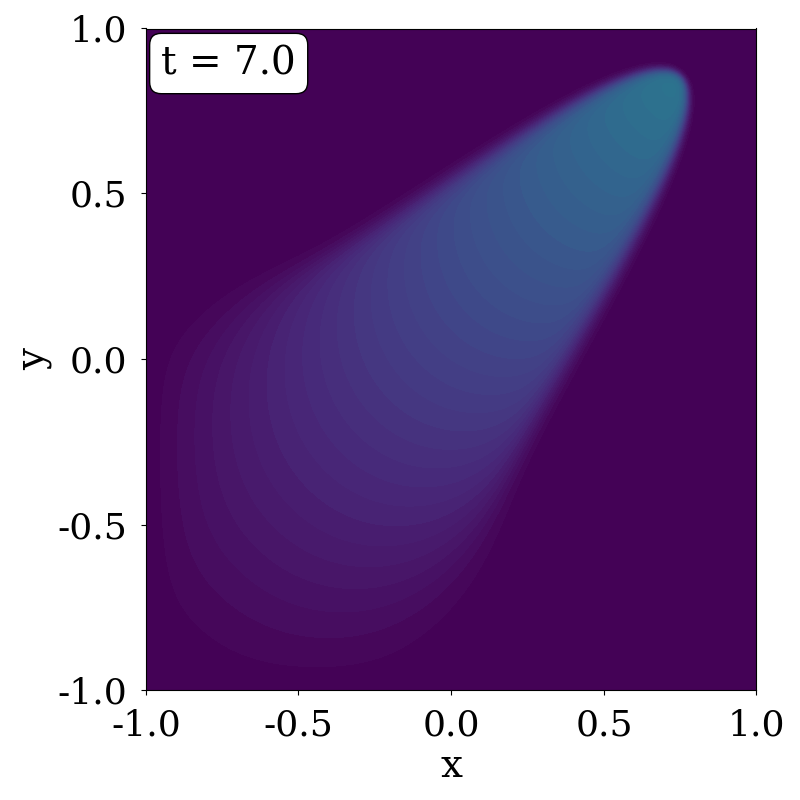}
            \subcaption{Burgers', $D = 3.25 \times 10^{-4}$}
        \end{minipage}
        \begin{minipage}{0.32\linewidth}
            \includegraphics[width=0.99\linewidth]{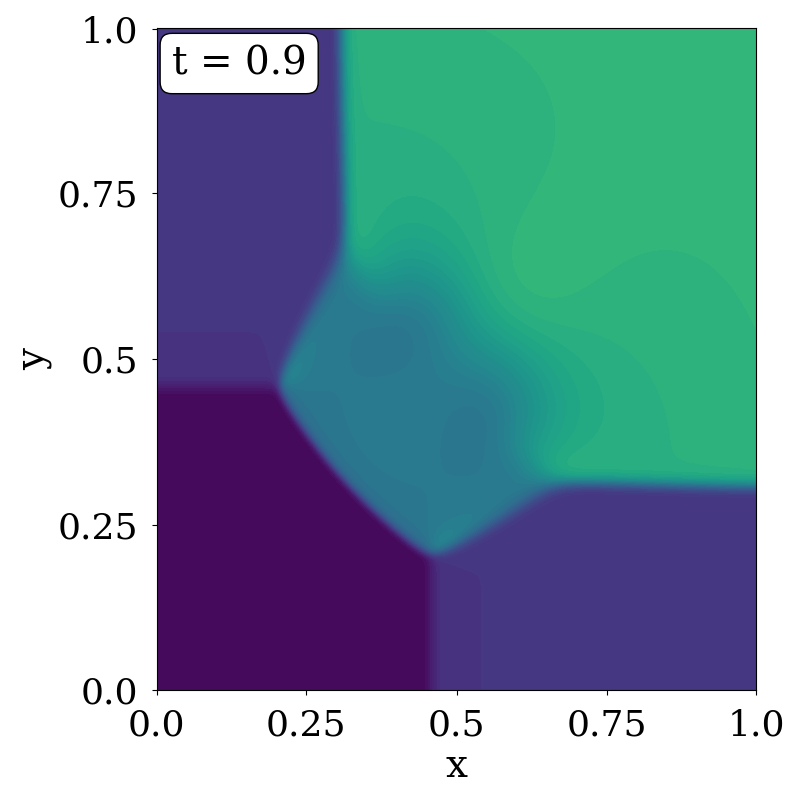}
            \subcaption{Euler, $p_4 = 1.0$}
        \end{minipage}
    }
    \caption{Sample solution snapshots at $t = 0$ (top),  $t = T/2$ (middle), and $t = T$ (bottom).}
    \label{fig:sample-snapshots}
\end{figure}

To assess the ability of the PROMs and the Schwarz alternating method in generating a robust model outside the training dataset, a single parameter is varied for each system and subsets of the resulting simulations are designated as the ``training'' dataset (from which all POD bases are computed) or ``testing'' dataset (from which predictive accuracy is measured). For the SWE, the Coriolis parameter $\mu$ is varied linearly from $-4.0$ to 0.0. For the Burgers' equations, the diffusion parameter $D$ is varied roughly logarithmically from 0.0001 to 0.001. For the Euler equations, the initial pressure in the upper right quadrant (noted here as $p_4$, though this corresponds to $p_1$ by the notation in~\cite{Schulz-Rinne1993}), is varied linearly from 0.5 to 1.5. The training and testing parameter values are summarized in Table~\ref{tab:param_settings}. Indicative snapshots of the resulting solution are shown in Figure~\ref{fig:sample-snapshots}. In all graphics measuring a quantity of interest across all parameter values, the training values are written in black text while testing values are written in red.

\begin{table}
    \caption{Training and testing parameters for each system.}
    \label{tab:param_settings}
    \centering
    \begin{tabular}{cccccccc}
    \toprule
    \multicolumn{2}{c}{2D SWE ($\mu$)} && \multicolumn{2}{c}{2D Burgers' ($D$)} && \multicolumn{2}{c}{2D Euler ($p_4$)} \\
    Train & Test && Train & Test && Train & Test \\
    \cmidrule{1-2}\cmidrule{4-5}\cmidrule{7-8}
    $-4.0$ & $-3.5$ && $1.00 \times 10^{-4}$ & $1.35 \times 10^{-4}$ && 0.50 & 0.625 \\
    $-3.0$ & $-2.5$ && $1.75 \times 10^{-4}$ & $2.50 \times 10^{-4}$ && 0.75 & 0.875 \\
    $-2.0$ & $-1.5$ && $3.25 \times 10^{-4}$ & $4.25 \times 10^{-4}$ && 1.00 & 1.125 \\
    $-1.0$ & $-0.5$ && $5.50 \times 10^{-4}$ & $7.50 \times 10^{-4}$ && 1.25 & 1.375 \\
    \phantom{$-$}0.0  &         && $1.00 \times 10^{-3}$ &          && 1.50 &       \\
    \bottomrule
    \end{tabular}
\end{table}

In measuring the accuracy of PROMs and coupled models, the error will be reported as the normalized space-time $\ell^2$ error computed for the $\varIdx$th variable as
\begin{equation}
    \epsilon_{\varIdx} \defEq \frac{\sum_{n=0}^{\numSnaps} \left\Vert \solVec^{\timeIdx} - \solVecRom^{\timeIdx} \right\Vert_2^2}{\sum_{n=0}^{\numSnaps} \left\Vert \solVec^{\timeIdx} \right\Vert_2^2}.
\end{equation}

\begin{figure}
    \centering
    \includegraphics[width=0.4\linewidth]{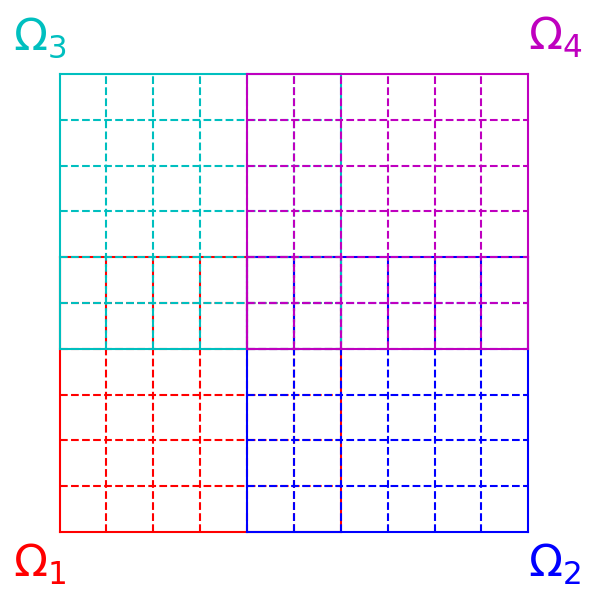}
    \caption{Exemplar $2 \times 2$ subdomain decomposition, with overlap regions of size $\numOverlap = 2$. Solid lines indicate subdomain boundaries, while dashed lines indicate mesh element edges. Different colors indicate different subdomains.}
    \label{fig:2x2_decomp}
\end{figure}

All coupled numerical experiments presented here are computed on a uniform $2\times2$ subdomain decomposition, with overlap regions of size $\numOverlap$ mesh elements between neighboring subdomains. This is roughly illustrated in Figure~\ref{fig:2x2_decomp}. This decomposition choice is the minimal decomposition to investigate coupling for waves propagating in the x- and y-directions; in general, finer decompositions will only increase runtime. For all experiments, the additive Schwarz alternating method is used to couple neighboring subdomains, and the absolute and relative Schwarz convergence tolerance is set to $\delta_{\text{abs}} = \delta_{\text{rel}} = 1 \times 10^{-11}$. Coupling is accomplished by a Dirichlet-Dirichlet interface between all neighboring subdomains, even for the ``non-overlapping'' ($\numOverlap = 0$) finite volume scenario discussed in Section~\ref{sec:schwarz_CCFV}. Additive Schwarz parallelism is implemented via {\tt OpenMP}, whereby each subdomain is tied to a single computational thread, and all share a common memory pool.

POD bases for monolithic and coupled PROMs are computed from the monolithic FOM solution. In the case of coupled PROMs, the monolithic FOM data is simply split according to the mesh decomposition. We note that monolithic FOM data is not always available (e.g., for multiphysics or assembly decompositions). In such a situation, the POD bases must necessarily be computed from coupled FOM solutions. In that case, datasets for each subdomain are naturally segregated and need not be manually split to compute POD bases. For the fluid flow systems investigated here, the decomposition is purely user-defined and the decomposed FOM solution is empirically observed to converge to the monolithic solution within machine precision, precluding the need to compute POD bases from decomposed FOM solutions. The centering vector $\solVecCent$ is computed as the snapshot average of the data matrix $\solMat$.

\subsection{PROM-PROM coupling}

\begin{figure}
    {
    \begin{minipage}{0.47\linewidth}
        \includegraphics[width=0.99\linewidth]{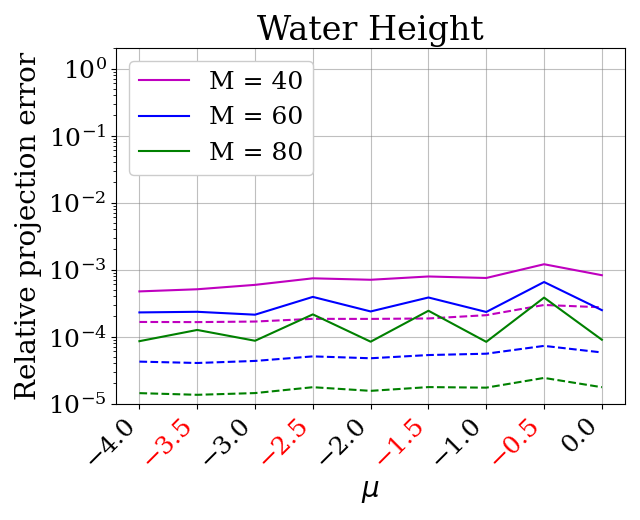}
    \end{minipage}
    \begin{minipage}{0.47\linewidth}
        \includegraphics[width=0.99\linewidth]{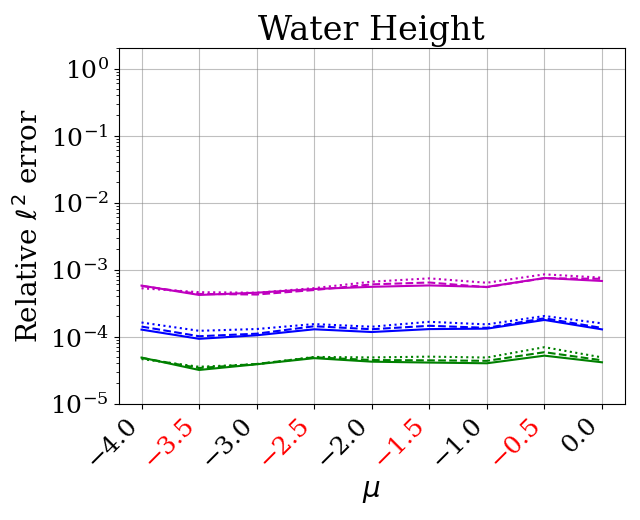}
    \end{minipage}
    \subcaption{SWE}
    }
    {
    \begin{minipage}{0.47\linewidth}
        \includegraphics[width=0.99\linewidth]{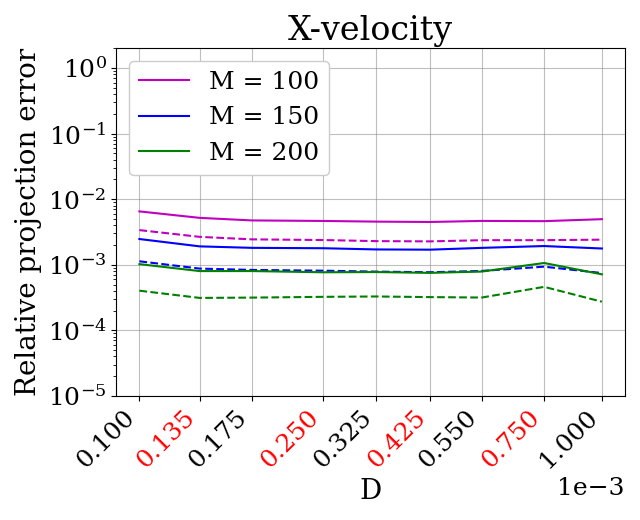}
    \end{minipage}
    \begin{minipage}{0.47\linewidth}
        \includegraphics[width=0.99\linewidth]{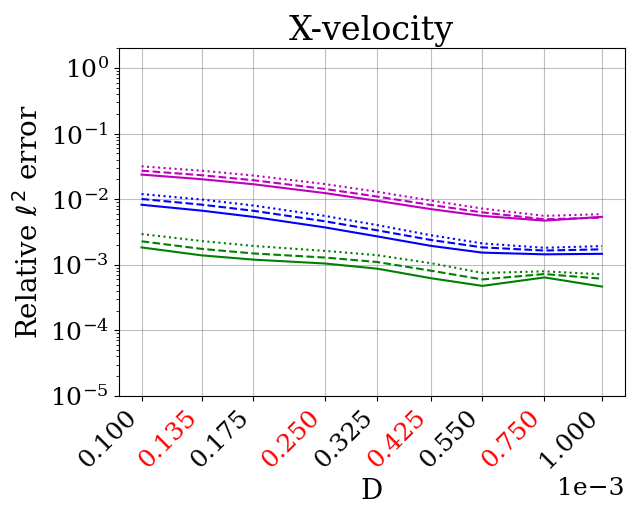}
    \end{minipage}
    \subcaption{Burgers'}
    }
    {
    \begin{minipage}{0.47\linewidth}
        \includegraphics[width=0.99\linewidth]{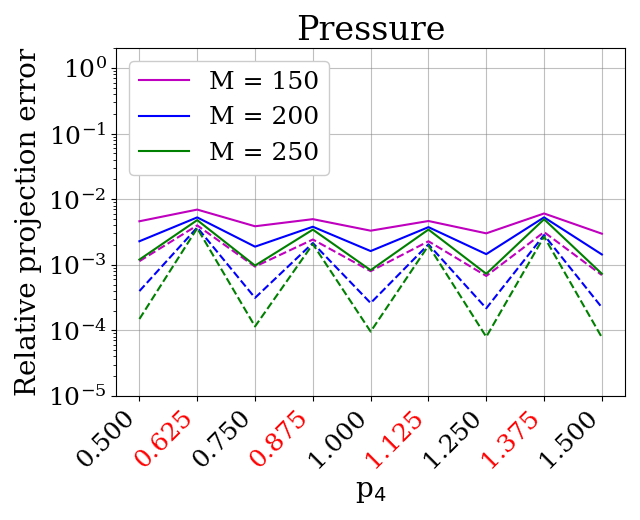}
    \end{minipage}
    \begin{minipage}{0.47\linewidth}
        \includegraphics[width=0.99\linewidth]{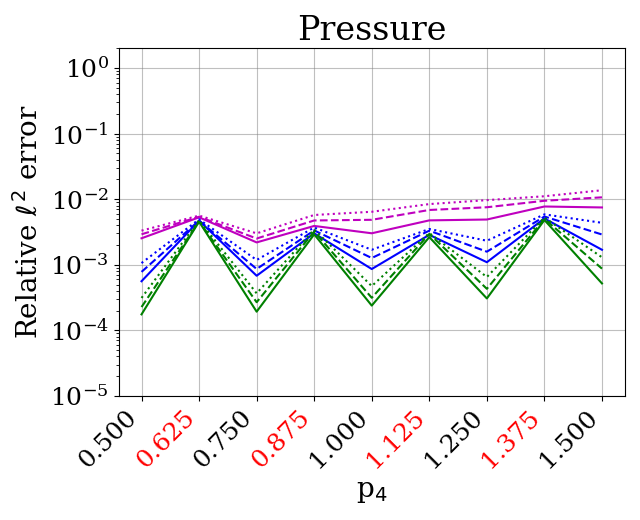}
    \end{minipage}
    \subcaption{Euler}
    \label{subfig:proj_error_euler}
    }
    \caption{Left column: relative space-time projection errors for monolithic solution (solid line) and decomposed solution with $\numOverlap = 0$ (dashed line), for fixed $\numModes$ in monolithic domain and each subdomain. Right column: online relative space-time error for coupled PROMs, for $\numOverlap = 0$ (solid line), $\numOverlap = 10$ (dashed line), and $\numOverlap = 20$ (dotted line). Training parameters are marked in black, while testing parameters are marked in red.}
    \label{fig:proj_and_unsampled_error}
\end{figure}

Before introducing hyper-reduction with the aim of achieving computational cost savings, it is important to investigate the characteristics of PROM coupling relative to monolithic PROM solutions. To begin, we examine the quality of the trial spaces induced by the monolithic and decomposed domains. The left column of Figure~\ref{fig:proj_and_unsampled_error} shows, for each problem, the projection error of the solution onto both the monolithic and decomposed trial spaces, for a fixed number of trial basis modes $\numModes$ in the monolithic domain and each subdomain. It is immediately clear that for a given basis size $\numModes$, the decomposed trial bases more accurately represent the solution as a whole, generally halving the error relative to the monolithic projection error at a given parameter value. This is an unsurprising result, as the subdomain-local trial basis may be optimized to better represent spatially-local dynamics. The decrease in error is smaller at higher basis sizes and test parameter values, as shown in the left image of Figure~\ref{subfig:proj_error_euler} for the Euler equations with $\numModes = 200$. This indicates the limits of the trial basis in generalizing to unseen parametric data for which the spatial locality is unable to compensate.

Of primary concern in constructing Schwarz-based solvers is the amount of overlap between neighboring subdomains. Classical Schwarz wisdom dictates that under sufficient convergence criteria, the Schwarz solution should converge to the monolithic solution regardless of the overlap size. However, a larger overlap generally results in faster convergence with fewer Schwarz iterations, at the cost of redundant solutions for mesh elements in the overlapping regions. In the case of coupling PROMs, one should never expect the decomposed PROM solution to converge exactly to the monolithic PROM solution, as the decomposed trial space is not the same as the monolithic trial space. One should expect, however, that the overlap size has relatively little effect on the solution accuracy. This is confirmed in the right column of Figure~\ref{fig:proj_and_unsampled_error}, which shows, for each problem, the online coupled PROM error as the overlap size and number of trial basis modes are varied. Indeed, the overlap size has a minimal effect on the solution accuracy when varied from zero overlap to 20 overlapping cells. In general, the zero overlap solution achieves the lowest error for any given parameter value, though this distinction is marginal.

\begin{figure}
    \begin{minipage}{0.32\linewidth}
        \includegraphics[width=0.99\linewidth]{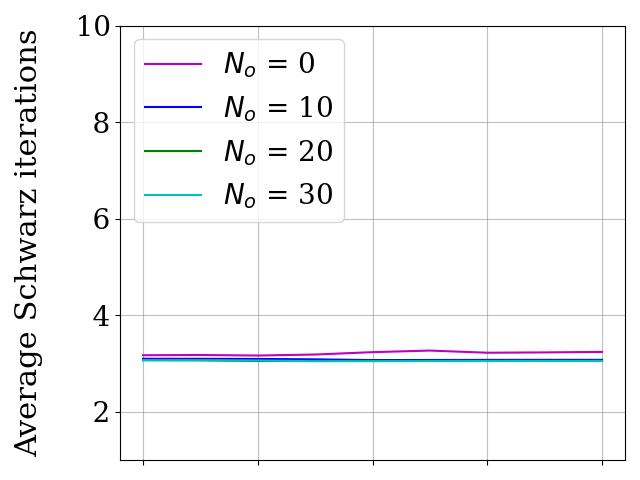}
    \end{minipage}
    \begin{minipage}{0.32\linewidth}
        \includegraphics[width=0.99\linewidth]{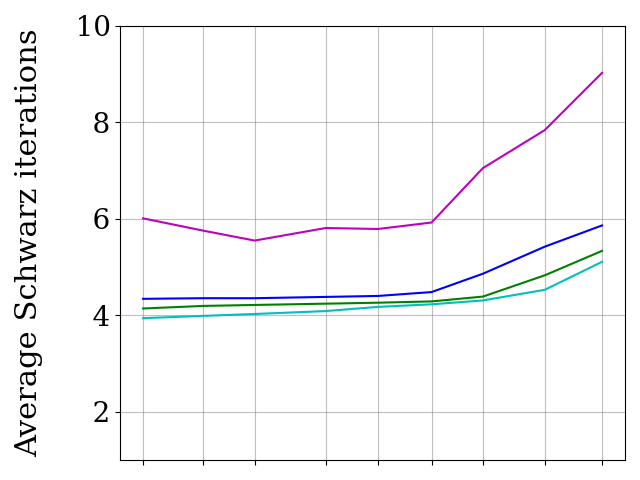}
    \end{minipage}
    \begin{minipage}{0.32\linewidth}
        \includegraphics[width=0.99\linewidth]{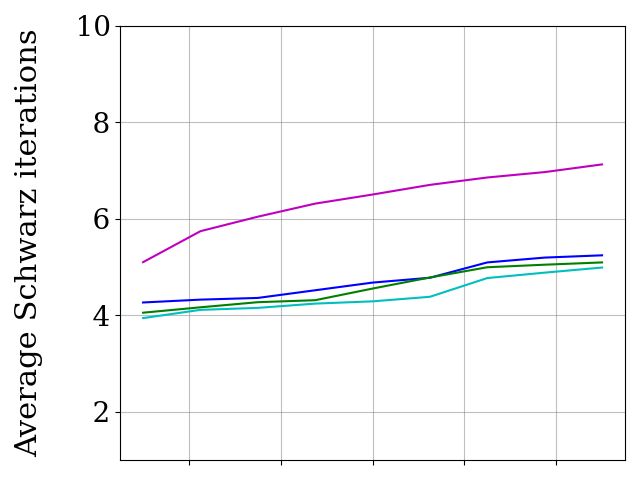}
    \end{minipage}

    \begin{minipage}{0.32\linewidth}
        \includegraphics[width=0.99\linewidth]{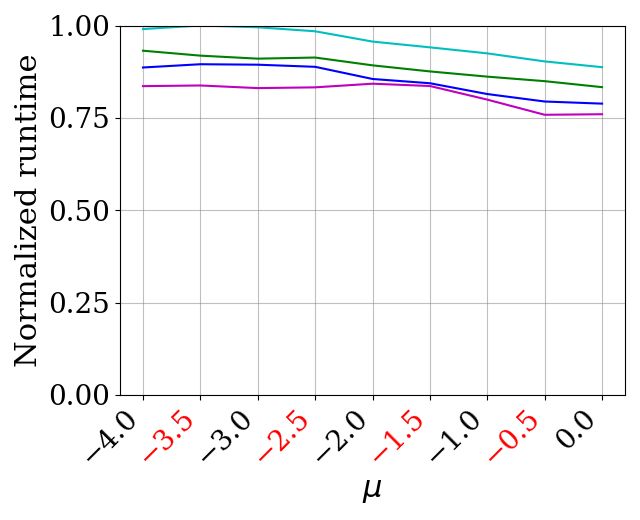}
        \subcaption{SWE, $\numModes = 80$}
        \label{subfig:schwarz_iters_swe}
    \end{minipage}
    \begin{minipage}{0.32\linewidth}
        \includegraphics[width=0.99\linewidth]{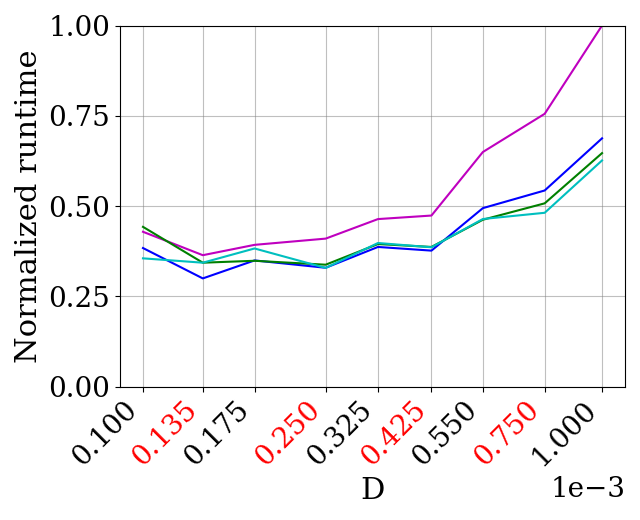}
        \subcaption{Burgers', $\numModes = 150$}
        \label{subfig:schwarz_iters_burgers}
    \end{minipage}
    \begin{minipage}{0.32\linewidth}
        \includegraphics[width=0.99\linewidth]{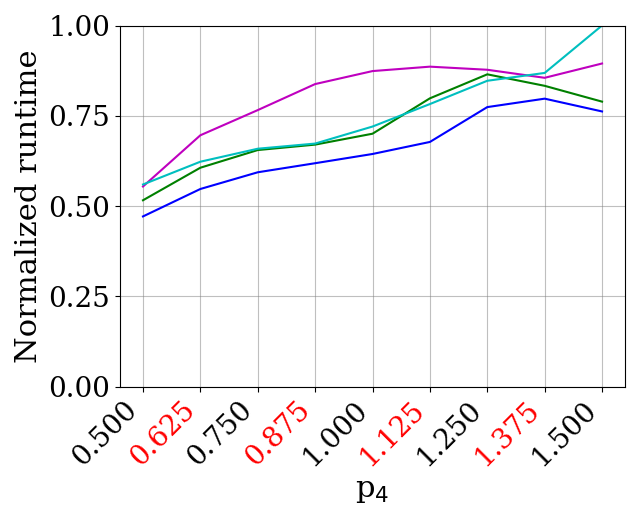}
        \subcaption{Euler, $\numModes = 200$}
        \label{subfig:schwarz_iters_euler}
    \end{minipage}
    
    \caption{Effect of overlap size on runtime performance for decomposed PROMs across a variety of overlap sizes $\numOverlap$. The top figures displays the average number of Schwarz iterations per time step, while the bottom figures show the normalized cost of each simulation. Training parameters are marked in black, while testing parameters are marked in red.}
    \label{fig:schwarz_iters}
\end{figure}

The effect of overlap size on Schwarz convergence speed is more significant, largely confirming the classical Schwarz wisdom, though this appears to be somewhat problem-dependent. The uppermost plots in Figure~\ref{fig:schwarz_iters} display the average number of Schwarz iterations for decomposed PROMs of varying overlap sizes and fixed trial basis sizes. For a fairly smooth problem such as the SWE, the overlap size has little effect, achieving convergence in roughly three iterations for every overlap size and parameter value. For Burgers' and the Euler equations, which are characterized by strong shocks propagating across Schwarz interfaces, Figures~\ref{subfig:schwarz_iters_burgers} and~\ref{subfig:schwarz_iters_euler} reveal not only a generally slower convergence (approximately four iterations for $\numOverlap > 0$), but particularly a sizable slowing of convergence for $\numOverlap = 0$. On average, this accounts for roughly two additional Schwarz iterations per time step. The effect is even stronger for higher diffusion coefficients $D$ (Burgers') and higher initial pressure $p_4$ in the upper right quadrant (Euler). 

It is interesting to remark that the increase in the average number of Schwarz iterations does not directly translate to a more expensive simulation. Increasing the overlap size naturally increases the computational cost for a single Schwarz iteration as each subdomain contains commensurately more mesh elements. The lower plots in Figure~\ref{fig:schwarz_iters} display normalized runtime measurements. For the SWE in Figure~\ref{subfig:schwarz_iters_swe}, the added cost of redundant cells in the overlap region outweighs any slim cost savings due to fewer Schwarz iterations, and the case where $\numOverlap = 0$ is both the most accurate and the fastest. For the Burgers' and Euler equations, the case where $\numOverlap = 0$ is usually (though not always) the most costly. While for the Burgers' equation at high $D$ this disparity is quite significant, for all other parameter values and for the Euler equations the costs are roughly comparable. On the whole, setting $\numOverlap = 0$ (the cell-centered finite volume Dirichlet-Dirichlet ``non-overlapping'' scenario) interestingly leads to an accurate coupled solution with minimal impact on computation expense, except in certain extremes. Thus, for the remainder of these results, we restrict analysis to the case for $\numOverlap = 0$.

\subsection{HPROM-HPROM coupling}

Recall that PROMs for sufficiently complex non-linear systems rarely achieve runtime reductions relative to the full-order model due the cost of evaluating the nonlinear terms in the FOM dimension $\numDOFs$. We now introduce hyper-reduction via collocation to achieve computational cost savings. Similar to the PROM results shown above, all coupled HPROM results are computed on a $2 \times 2$ decomposed domain and utilize an HPROM model in every subdomain. For all HPROM domains (both monolithic and decomposed), the trial basis size is fixed at a sufficient level to at least guarantee low projection error. The SWE experiments are fixed at $\numModes = 80$, the Burgers' equation at $\numModes = 150$, and the Euler equations at $\numModes = 200$.

In generating the sample meshes, the sample indices set $\mathcal{S}$ is first seeded with any Schwarz interface samples as dictated by the interface sampling interval $\numBound$. The Burgers' equation sample meshes are further seeded with the degrees of freedom selected by QDEIM~\cite{Drmac2016}, computed from a trial basis of size $\numModes = 200$. The remainder of the sample points under a fixed sample mesh size $\numSamps$ are selected randomly, with higher sampling rates simply appending additional sampled mesh elements to those selected at lower sampling rates, ensuring uniformity in all numerical experiments. A range of sample mesh sizes are investigated, ranging roughly logarithmically from 0.1\% to 10\% of the full mesh size, assuming the number of seed points does not already exceed this sample mesh size. In general, all HPROM solutions for the SWE are stable above $\numSamps = 0.1\% \times \numDOFs$, while solutions for the Burgers' and Euler equations requires a sampling rate of more than $1.0\% \times \numDOFs$. This is a natural consequence of the system dynamics. Whereas the SWE generate a fairly smooth solution which is easily represented by a sparse sample mesh, the Burgers' and Euler equations require relatively dense sample meshes to accurately propagate shocks.

\begin{figure}
    \begin{minipage}{0.32\linewidth}
        \includegraphics[width=0.99\linewidth]{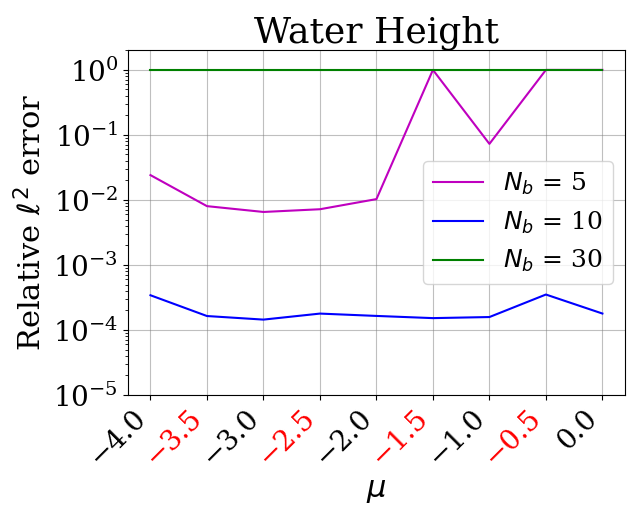}
        \subcaption{SWE, $\numModes = 80$, $\numSamps = 0.5\% \times \numDOFs$}
    \end{minipage}
    \begin{minipage}{0.32\linewidth}
        \includegraphics[width=0.99\linewidth]{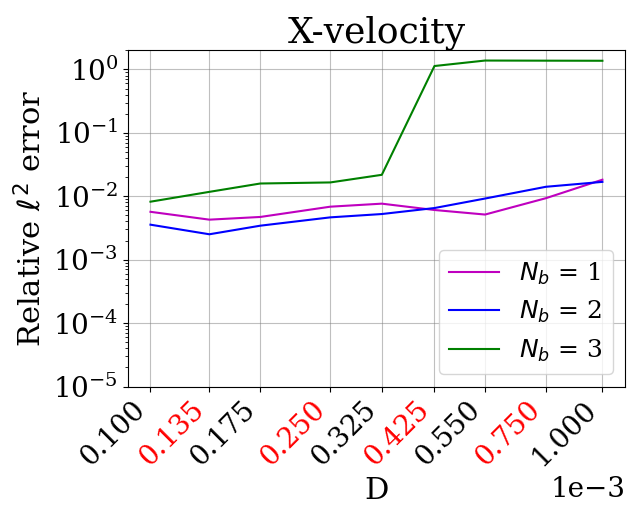}
        \subcaption{Burgers', $\numModes = 150$, $\numSamps = 3.75\% \times \numDOFs$}
    \end{minipage}
    \begin{minipage}{0.32\linewidth}
        \includegraphics[width=0.99\linewidth]{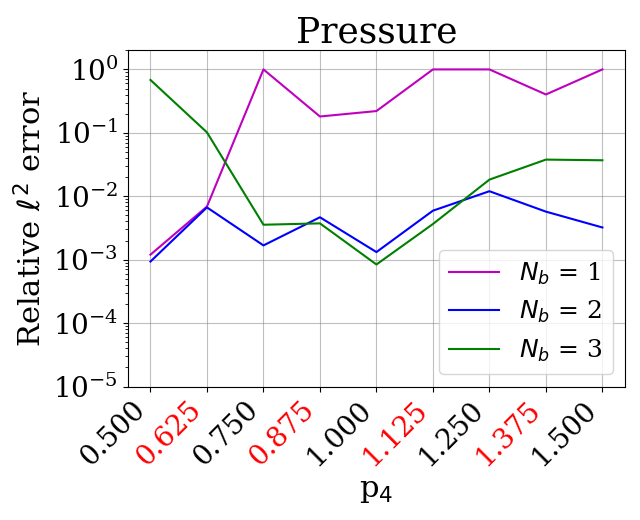}
        \subcaption{Euler, $\numModes = 200$, $\numSamps = 5.0\% \times \numDOFs$}
    \end{minipage}
    \caption{Online relative space-time error of HPROMs for various interface sampling intervals $\numBound$, with fixed $\numModes$ and $\numSamps$. Training parameters are marked in black, while testing parameters are marked in red.}
    \label{fig:hrom_error_overlap}
\end{figure}

The body of hyper-reduction literature confirms the intuition that increasing the sample mesh size generally improves solution accuracy and robustness. For coupled HPROMs, however, the question of interface sampling is more nuanced. For all cases investigated here, failing to seed the sample mesh with interface mesh elements leads to an unstable solution, except at extremely high sampling random rates (10\% $\times \numDOFs$ for the SWE). We must then deliberately seed the sample mesh with elements at the Schwarz interface, in this case utilizing a simple fixed sampling interval $\numBound$, as described in Section~\ref{subsec:bound_samp}. As stated previously, decreasing the sampling interval naturally improves the ease of transmitting interface information between neighboring subdomains, but draws sample mesh elements away from the subdomain interior and degrades overall solution accuracy. This balancing act is illustrated in Figure~\ref{fig:hrom_error_overlap}, where an optimal interface sampling interval generates accurate and stable coupled HPROMs, while values of $\numBound$ with are either too high or too low may generate unstable solutions. For the SWE, this optimum appears to be around $\numBound = 10$, while for the Burgers' and Euler equations this is approximately $\numBound = 2$, or every other interface cell. The disparity between these optimums is again due to the smoothness of the SWE solution, in contrast to the shocks experienced by the Burgers' and Euler solutions. The ability to accurately transmit a shock across the Schwarz interface reasonably requires a much higher sample mesh resolution.

\begin{figure}
    \begin{subfigure}{0.99\linewidth}
        \includegraphics[width=0.99\linewidth]{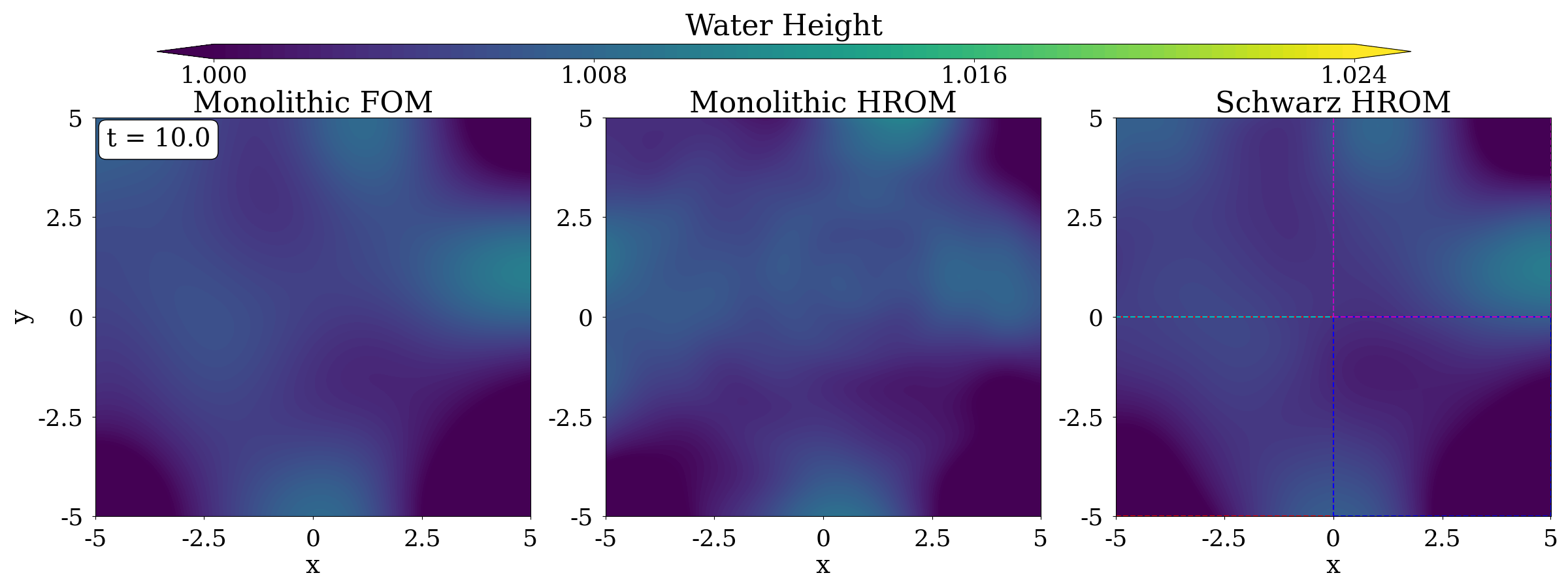}
        \subcaption{SWE, $\mu = -0.5$, $\numModes = 80$, $\numSamps = 0.5\% \times \numDOFs$, $\numBound = 10$}
    \end{subfigure}
    \begin{subfigure}{0.99\linewidth}
        \includegraphics[width=0.99\linewidth]{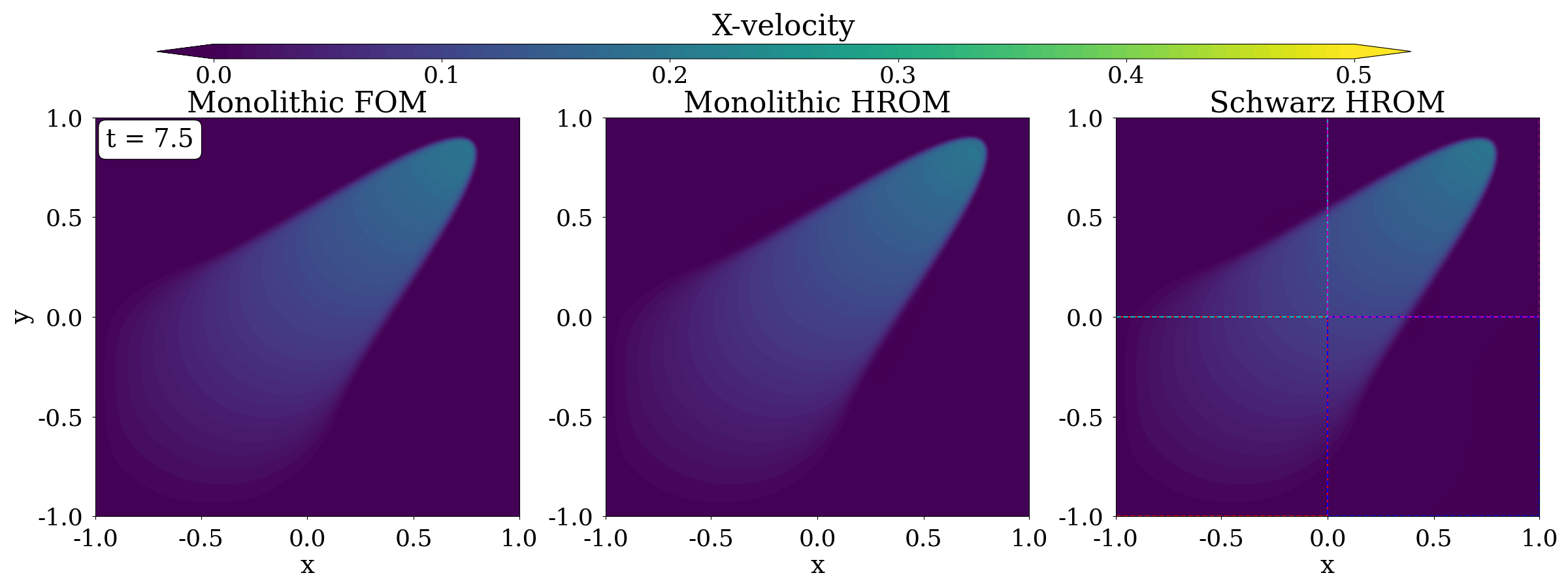}
        \subcaption{Burgers', $D = 1.35 \times 10^{-4}$, $\numModes = 150$, $\numSamps = 3.75\% \times \numDOFs$, $\numBound = 2$}
    \end{subfigure}
    \begin{subfigure}{0.99\linewidth}
        \includegraphics[width=0.99\linewidth]{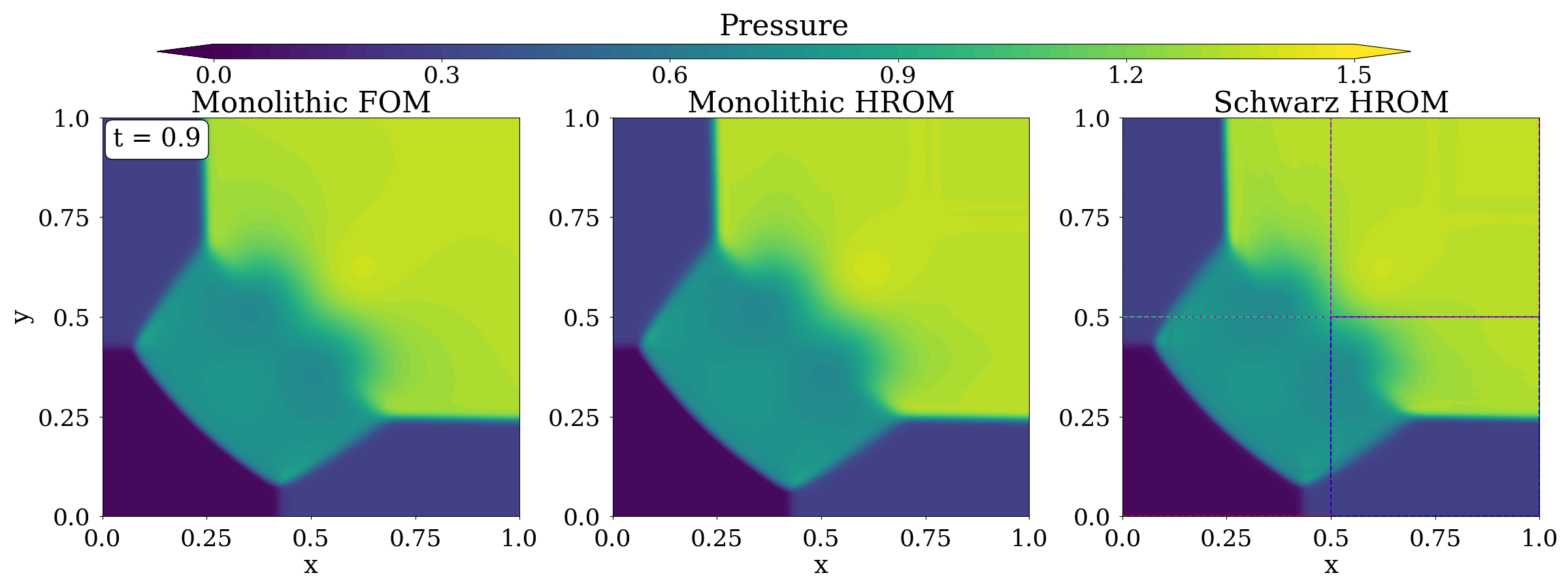}
        \subcaption{Euler, $p_4 = 1.375$, $\numModes = 200$, $\numSamps = 5.0\% \times \numDOFs$, $\numBound = 2$}
    \end{subfigure}
    \caption{Monolithic FOM (left), monolithic HPROM (center), and coupled HROM with $\numOverlap = 0$ (right) solution contours at $t = T$, with fixed $\numModes$ and $\numSamps$. Subdomain interfaces are marked by dashed lines.} 
    \label{fig:hrom_contours}
\end{figure}

Under sufficient sampling rates for a stable solution, it is useful to know whether the Schwarz coupling \textit{induces} errors at interfaces. That is, as each subdomain has a unique trial space which does not guarantee perfect continuity across interfaces, are there inordinately high error levels at the interfaces? Or are the Schwarz convergence criteria sufficient to ensure accurate transmission between disparate trial spaces? For a fixed trial basis dimension $\numModes$ and sample mesh size $\numSamps$, Figure~\ref{fig:hrom_contours} displays representative monolithic FOM, monolithic HPROM, and decomposed HPROM solutions at $t = T$. In a broad sense, the decomposed HPROM solution mirrors that of the monolithic HPROM closely: slight ``ringing'' effects are apparent in the Burgers' solution, and more so in the Euler solution. The decomposed SWE solution roughly appears to approximate the FOM solution more closely than does the monolithic HPROM. Most notably, there are not any apparent discontinuities at the subdomain boundaries which would signal poor convergence at the interface. 

\begin{figure}
    \centering
    \begin{minipage}{0.32\linewidth}
        \includegraphics[width=0.99\linewidth]{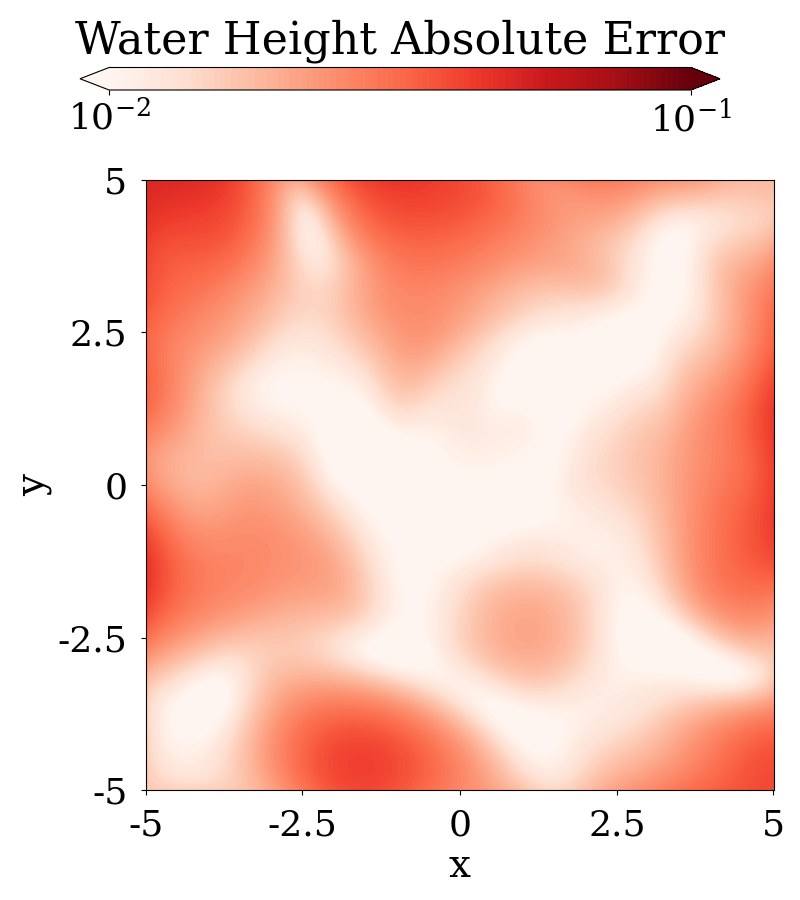}
    \end{minipage}
    \begin{minipage}{0.32\linewidth}
        \includegraphics[width=0.99\linewidth]{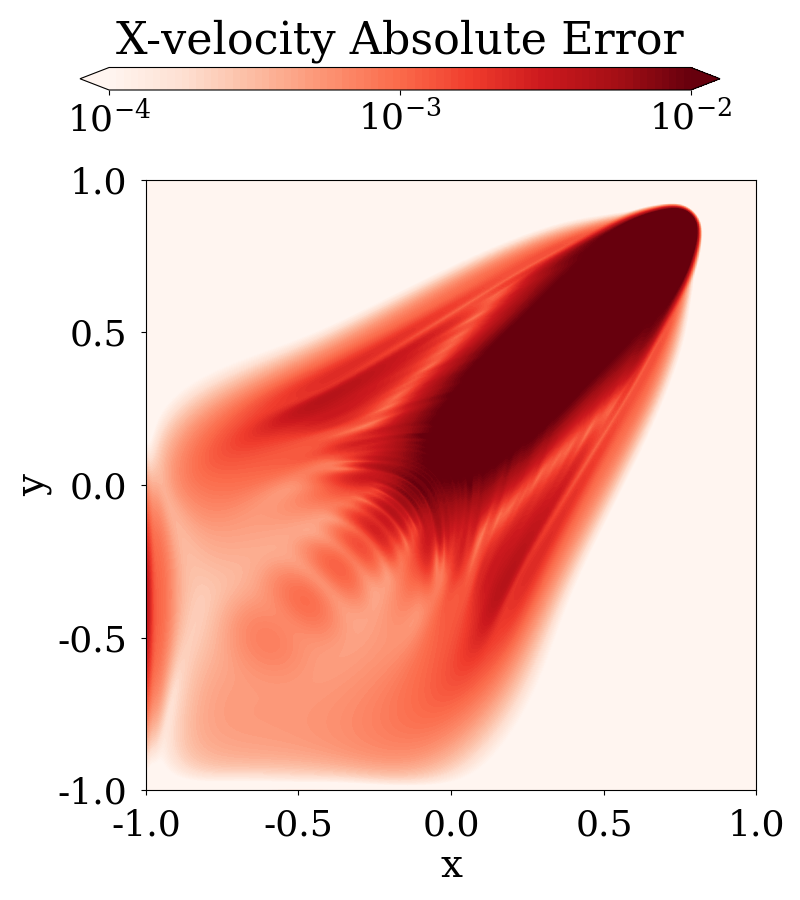}
    \end{minipage}
    \begin{minipage}{0.32\linewidth}
        \includegraphics[width=0.99\linewidth]{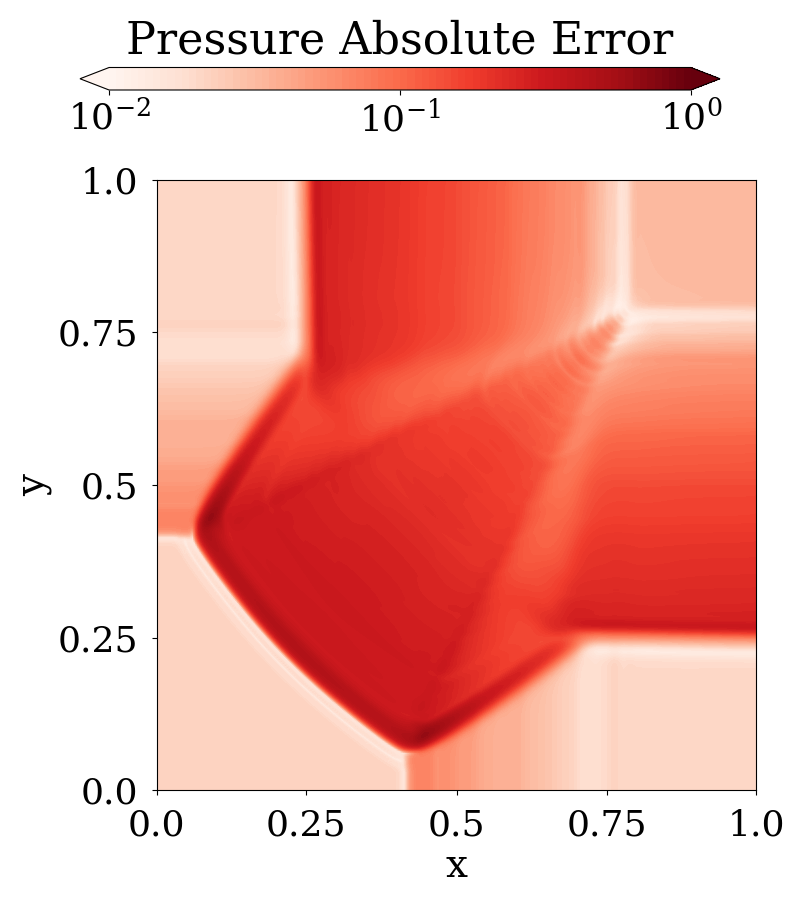}
    \end{minipage}

    \begin{minipage}{0.32\linewidth}
        \includegraphics[width=0.99\linewidth]{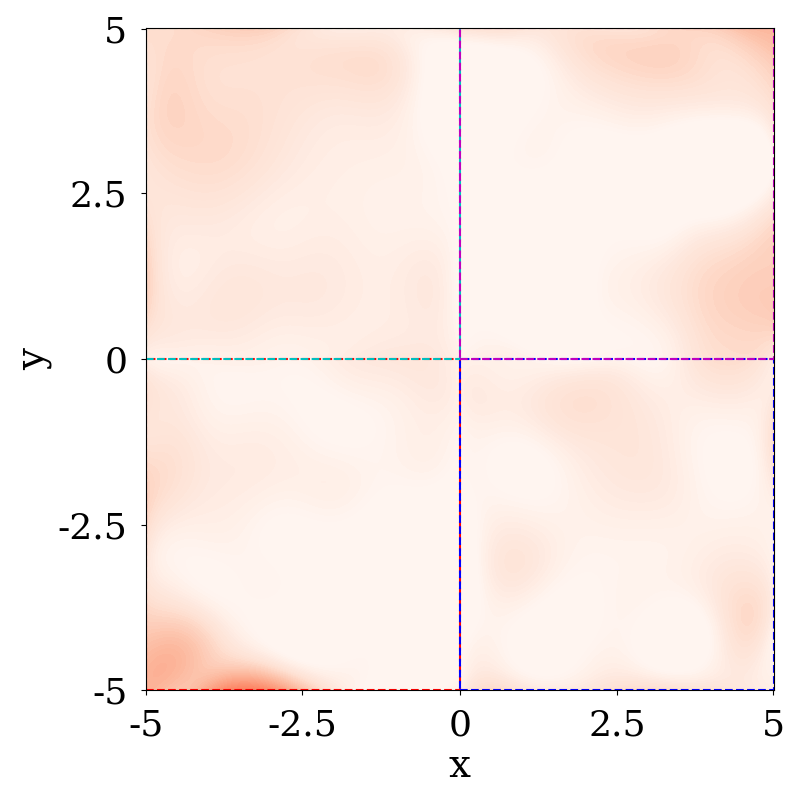}
        \subcaption{\centering SWE, $\mu = -0.5$, \newline $\numModes = 80$, $\numSamps = 0.5\% \times \numDOFs$, $\numBound = 10$}
    \end{minipage}
    \begin{minipage}{0.32\linewidth}
        \includegraphics[width=0.99\linewidth]{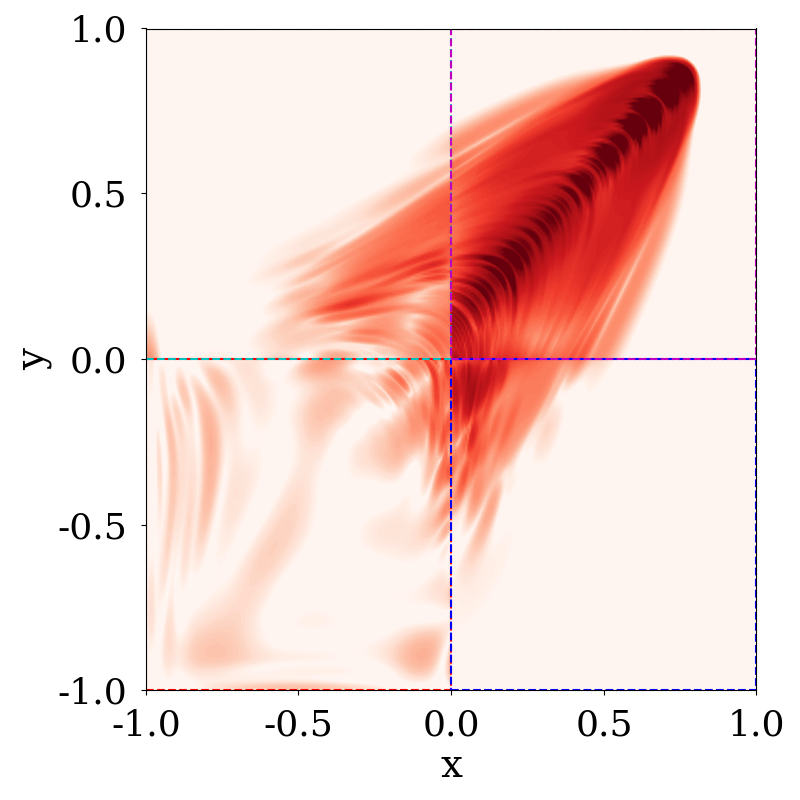}
        \subcaption{\centering Burgers', $D = 1.35 \times 10^{-4}$,\newline $\numModes = 150$, $\numSamps = 3.75\% \times \numDOFs$, $\numBound = 2$}
    \end{minipage}
    \begin{minipage}{0.32\linewidth}
        \includegraphics[width=0.99\linewidth]{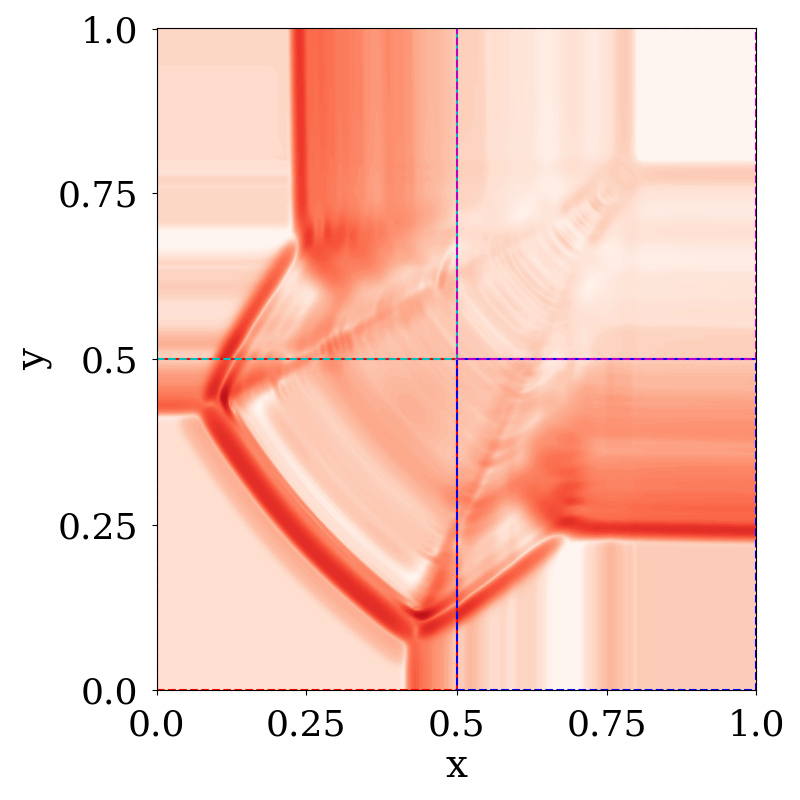}
        \subcaption{\centering Euler, $p_4 = 1.375$, \newline $\numModes = 200$, $\numSamps = 5.0\% \times \numDOFs$, $\numBound = 2$}
    \end{minipage}
    
    \caption{Average absolute spatial error fields for representative monolithic (top) and decomposed (bottom) HPROM solutions with $\numOverlap = 0$. Subdomain interfaces are marked by dashed lines. Note that the computational cost of the monolithic and decomposed solutions are not equivalent; see Figure~\ref{fig:pareto} for comparisons.} 
    \label{fig:space_avg_error}
\end{figure}

To examine the error behavior more precisely, Figure~\ref{fig:space_avg_error} shows the time-average absolute error for the same monolithic and decomposed HPROM solutions. If the Schwarz coupling induces error at the interface, one might expect high error concentration along the subdomain boundaries. This does not appear to be the case, as not only are errors generally lower for the decomposed solution (again, the greater accuracy of subdomain-local trial bases), but error \textit{patterns} are largely the same (dominated by features of the flow dynamics, such as shocks). This indicates that the Schwarz coupling algorithm has no difficulty in accurately transmitting information under sufficient interface sampling conditions.

\begin{figure}
    \begin{minipage}{0.49\linewidth}
        \includegraphics[width=0.99\linewidth]{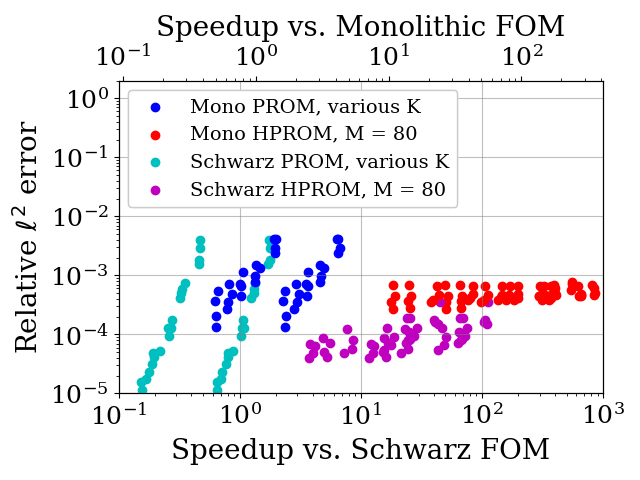}
        \subcaption{SWE, $\numBound = 10$}
    \end{minipage}
    \begin{minipage}{0.49\linewidth}
        \includegraphics[width=0.99\linewidth]{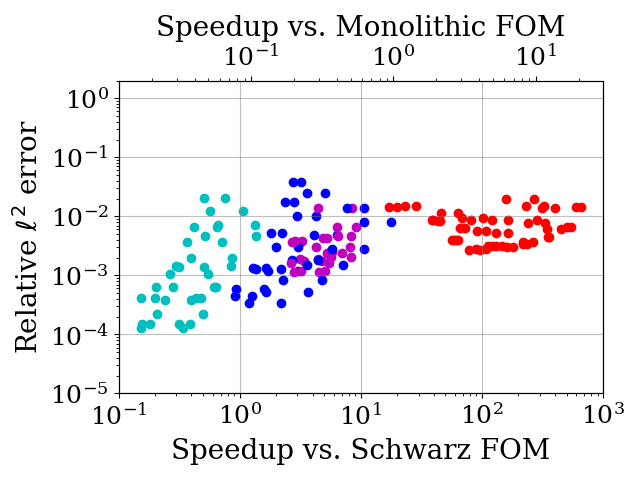}
        \subcaption{Burgers', $\numBound = 2$}
    \end{minipage}
    \begin{center}
    \begin{minipage}{0.49\linewidth}
        \includegraphics[width=0.99\linewidth]{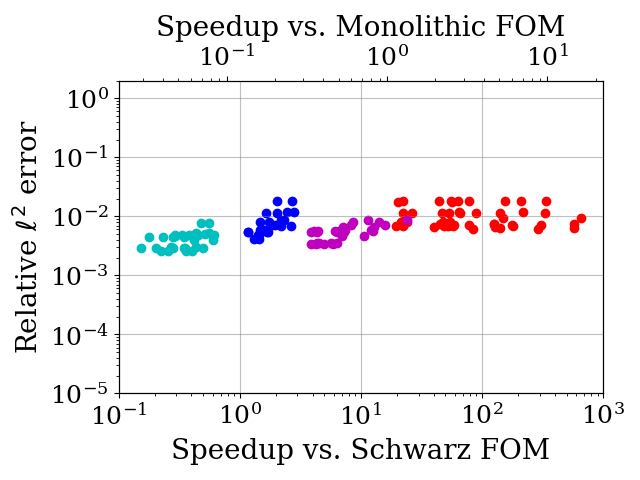}
        \subcaption{Euler, $\numBound = 2$}
    \end{minipage}
    \end{center}
    \caption{Pareto fronts comparing trade-offs in computational accuracy and cost reduction (lower-right is better) under a variety of simulation configurations, for all \textit{test} parameter values. All decomposed solutions are computed with $\numOverlap = 0$. The top $x$-axis measures speedup relative to the monolithic full-order solution, while the bottom $x$-axis is relative to the decomposed all-FOM solution with $\numOverlap = 0$.}
    \label{fig:pareto}
\end{figure}

Finally, we address the effect of coupling on the cost efficiency of HPROMs. Figure~\ref{fig:pareto} plots the solution error as a function of either the computational speedup with respect to the all-FOM decomposed solution or the monolithic FOM solution. For the PROMs, each point represents a different test parameter value and trial basis size $\numModes$. For the HPROMs, each point represents a different test parameter value and sample mesh size $\numSamps$. These cost measurements take into account the three additional cores used to parallelize the additive Schwarz algorithm.

Unsurprisingly, the PROMs are generally unable to achieve any computational cost savings compared to the monolithic PROM solution, and the decomposed PROM solution is invariably slower than the monolithic PROM solution due to the added cost of Schwarz iterations. The HPROM solutions, on the other hand, are all capable of achieving noticeable speedups relative to the decomposed FOM solution. The decomposed HPROM solutions broadly achieve lower error levels than the monolithic HPROM solutions, but do not reach the same maximum speedup levels. For the SWE, the monolithic HPROM is able to achieve speedups over 1,000 times faster than the decomposed FOM, while the decomposed HPROM is able to achieve roughly 100 times speedup at best. The same is true for the Burgers' and Euler equations, where the monolithic and decomposed HPROM solutions achieve a maximum speedup of 100 times and 10 times, respectively, over the decomposed FOM. Note that especially for the shock-driven problems, the decomposed HPROMs are not able to achieve significant speedups over the monolithic FOM, as they require relatively dense sample meshes to capture the shock dynamics and low interface sampling intervals to transmit shocks across subdomain boundaries. That being said, these results indicate that the Schwarz alternating method can robustly couple HPROMs and achieve good cost savings relative to an equivalent coupled FOM solution. For systems which \textit{must} be coupled, Schwarz coupling may offer appreciable computational speedups.

% Conclusions
\section{Conclusions} \label{sec:conc}

In this work, the Schwarz alternating method has been presented as a
minimally-intrusive method for coupling data-driven projection-based reduced-order models. By the iterative application of appropriate interface conditions, PROMs may be applied to arbitrary decompositions of a spatial domain and seamlessly joined to achieve accurate and lower cost solutions for dynamical systems. This has been demonstrated for three nonlinear 2D hyperbolic fluid flow problems, all run in the predictive regime, for which the transmission of propagating waves and representation of shocks is extremely challenging. These numerical experiments reinforce the fact that subdomain-local trial spaces are capable of more accurately representing complex flow phenomena with fewer basis modes than equivalent monolithic solutions. With an appropriate subdomain interface sampling strategy, hyper-reduced PROMs are also capable of delivering stable and accurate solutions which achieve up to two orders of magnitude speedup over equivalent coupled FOM solutions. Although speedups relative to a monolithic FOM are more difficult to achieve, the method still stands to accelerate applications which must be coupled, such as mechanical assemblies or multiphysics systems. Further, under the specific discretization choice of the cell-centered finite volume method, these results also indicate that a ``non-overlapping'' Dirichlet-Dirichlet Schwarz coupling where no physical cells overlap is capable of achieving convergence. This contrasts with traditional vertex-centered Dirichlet-Neumann non-overlapping interfaces, and hints at interesting possibilities in simplified discretization and coupling.

There are several avenues of future research which may follow from this work. Experiments with non-conformal meshes and variable time steps may further demonstrate the versatility of the Schwarz alternating method. Automated selection of the appropriate model (e.g., trial basis size, sample mesh) for each subdomain may further improve accuracy and computational speedups. In particular, the time adaptation of interface samples may alleviate many of the issues observed for the shock-driven Burgers' and Euler systems investigated here. It may be possible to speed up the convergence of this Schwarz-based approach through the development of optimized transmission BCs. Finally, the incorporation of alternative modeling strategies (such as non-intrusive ROMs) may bring additional modeling flexibility and cost savings.

\section*{Acknowledgements} \label{sec:acknowl}

Support for this work was received through Sandia National Laboratories' Laboratory Directed Research and Development (LDRD) program and through the U.S. Department of Energy, Office of Science, Office of Advanced Scientific Computing Research, Mathematical Multifaceted Integrated Capability Centers (MMICCs) program, under Field Work Proposal 22025291 and the Multifaceted Mathematics for Predictive Digital Twins (M2dt) project. Additionally, the writing of this manuscript was funded in part by the fourth author’s (Irina Tezaur’s) Presidential Early Career Award for Scientists and Engineers (PECASE).

Sandia National Laboratories is a multi-mission laboratory managed and operated by National Technology and Engineering Solutions of Sandia, LLC., a wholly owned subsidiary of Honeywell International, Inc., for the U.S. Department of Energy’s National Nuclear Security Administration under contract DE-NA0003525.

The authors wish to thank Alejandro Mota and Daria Koliesnikova for engaging in enlightening discussions related to details pertaining to the Schwarz alternating method for domain decomposition-based coupling.

% Appendices...
%\input{AppendixA}

\printbibliography

\end{document}